\documentclass[opre,nonblindrev]{informs4}
\usepackage{apacite}

\OneAndAHalfSpacedXI 

\usepackage{endnotes}
\let\footnote=\endnote

\usepackage{natbib}
 \bibpunct[, ]{(}{)}{,}{a}{}{,}%
 \def\bibfont{\small}%
\TheoremsNumberedThrough     
\ECRepeatTheorems

\EquationsNumberedThrough    

\usepackage{bigstrut,bigdelim,multirow}
\usepackage{booktabs}
\usepackage{enumerate}
\usepackage{url}
\usepackage{verbatim}
\usepackage[ruled,vlined,linesnumbered]{algorithm2e}
\SetKwInput{KwInput}{Input}
\SetKwInput{KwOutput}{Output} 
\usepackage{amsmath}
\usepackage{bm}
\usepackage{lmodern}
\usepackage{xcolor}

\usepackage{appendix}

\usepackage{graphicx}
\usepackage{subcaption}
\usepackage{float} 

\usepackage[colorlinks,linkcolor=blue,anchorcolor=blue,citecolor=blue]{hyperref}
\begin{document}

\graphicspath{{figures/}}

\def\tinyl{\mbox{\tiny L}}
\def\tinyf{\mbox{\tiny F}}
\def\tinyu{\mbox{\tiny U}}




\TITLE{Optimization with Multi-sourced Information and Unknown Reliability: A Distributionally Robust Approach}


\ARTICLEAUTHORS{%
\AUTHOR{Yanru Guo}
\AFF{Department of Industrial and Operations Engineering, University of Michigan, Ann Arbor, MI 48109, \EMAIL{yanruguo@umich.edu}} 
\AUTHOR{Ruiwei Jiang}
\AFF{Department of Industrial and Operations Engineering, University of Michigan, Ann Arbor, MI 48109, \EMAIL{ruiwei@umich.edu}} 
\AUTHOR{Siqian Shen}
\AFF{Corresponding author; Department of Industrial and Operations Engineering, University of Michigan, Ann Arbor, MI 48109, \EMAIL{siqian@umich.edu}}
} 

\ABSTRACT{%
In problems that involve input parameter information gathered from multiple data sources with varying reliability, incorporating decision makers' trust on different sources in optimization models can potentially improve solution performance. In this work, we propose a novel multi-reference distributionally robust optimization (MR-DRO) framework, where the model inputs are uncertain and their probability distributions can be statistically inferred from multiple information sources. Via nonparametric data fusion, we construct a Wasserstein ambiguity set to minimize the worst-case expected cost of a stochastic objective function, accounting for both uncertainty and unknown reliability of several given information sources. We reformulate the MR-DRO model as a linear program given linear objective and constraints in the original problem. We also incorporate a dynamic trust update mechanism that adjusts the trust for each source based on its performance over time. In addition, we introduce the concept of probability dominance to identify sources with dominant trust. Via computational studies using resource allocation and portfolio optimization instances, we demonstrate the effectiveness of the MR-DRO approach compared to traditional optimization frameworks relying on a single data source. Our results highlight the significance of integrating (dynamic) decision maker's trust in optimization under uncertainty, particularly when given diverse and potentially conflicting input data. 
}%

\KEYWORDS{(R) Robust Optimization; multi-sourced stochastic optimization; nonparametric distributionally robust optimization; convex programming; probabilistic dominance}

\maketitle

%


\newpage
\section{Introduction}\label{sec:intro}

Information uncertainty is ubiquitous in real-world decision-making problems. Such uncertainty often arises from limited data, measurement noise, and prediction errors generated by multiple information sources. Classical stochastic and robust optimization frameworks typically assume that uncertainty can be characterized a priori either by a known probability distribution or by a predefined uncertainty set constructed from available data. These approaches generally do not differentiate among data sources in terms of reliability or account for how heterogeneous data quality affects solution robustness.

In stochastic optimization \citep{birge2011introduction}, the true distribution of uncertain parameters is assumed to be known, and decisions are optimized with respect to expected performance. Robust optimization \citep{ben2002robust, bertsimas2011theory}, by contrast, optimizes against the worst-case realization within a prescribed uncertainty set. Distributionally robust optimization (DRO) \citep{kuhn2025distributionally} bridges these paradigms by constructing an ambiguity set of probability distributions consistent with observed data and optimizing against the worst-case distribution within that set \citep{calafiore2006distributionally, delage2010distributionally, wiesemann2014distributionally, kuhn2024distributionallyrobustoptimization, jiang2016data}. Among various constructions, Wasserstein-metric-based ambiguity sets \citep{KR:58} have become particularly prominent due to their statistical convergence guarantees and computational tractability \citep{mohajerin2018data, xie2021distributionally, shen2023chance, chen2024data}. In most existing work, however, the Wasserstein ball is centered at a single empirical distribution derived from one data source, implicitly assuming homogeneous data quality and limiting applicability in settings with multiple heterogeneous sources.

In practice, decision-makers often rely on diverse information sources with varying levels of accuracy and bias. Weather forecasting integrates satellite, radar, and ground-station measurements; healthcare decisions synthesize signals from different diagnostic devices; and cyber-physical systems aggregate distributed sensor data. While data fusion techniques have been extensively studied in engineering and AI \citep{liggins2017handbook, albahri2023systematic, ounoughi2023data}, their integration into optimization-based decision models under exogenous uncertainty remains limited in the OR/MS literature.

Recent DRO-based approaches begin to incorporate multi-source information. \citet{rychener2024} construct ambiguity sets as intersections of Wasserstein balls centered at source-specific empirical distributions to hedge against source bias. Although this framework mitigates misspecification at the source level, it assumes static source distributions and sufficient proximity to the target distribution—assumptions that may fail under distribution shifts or substantial cross-source heterogeneity. \citet{guo-etal-cdc-2024} propose a parametric data fusion method that aggregates trust-weighted moments to form a Normal-based ambiguity set, yielding a tractable linear reformulation under the 1-Wasserstein metric. However, the reliance on a common parametric family restricts applicability when source distributions differ structurally. More recently, \citet{keehan2025don} develop a Wasserstein DRO framework for temporally evolving multi-source data using weighted empirical distributions and concentration-based weight calibration. While this advances dynamic integration, it relies on bounded Wasserstein drift assumptions and requires nontrivial joint calibration of weights and ambiguity radii, particularly in high-dimensional settings.

These limitations reveal a fundamental gap: there is a need for optimization frameworks that (i) flexibly integrate heterogeneous data sources without restrictive parametric assumptions, (ii) explicitly model varying and evolving source reliability, and (iii) maintain tractability and robust performance guarantees.

To address this gap, we introduce a multi-reference distributionally robust optimization (MR-DRO) framework that incorporates trust as an endogenous measure of source reliability. The concept of trust, widely studied in human–robot interaction \citep{lee2004trust, lewis2018role, 6005228, 5753439}, captures a decision maker’s belief in the reliability of each information source. In our setting, trust governs how source data are fused into a nonparametric empirical reference distribution. Unlike existing parametric approaches, our framework does not require sources to share a common distributional family, nor does it assume proximity to a true underlying distribution.

We further develop a dynamic trust update mechanism that adjusts source weights over time based on realized performance. Drawing inspiration from Bayesian trust learning in human–robot interaction \citep{guo2021modeling, guo2021reverse} and adaptive expert-weighting methods in online learning \citep{herbster1998tracking}, our mechanism evaluates alignment between historical source data and observed realizations to recalibrate trust levels. This enables the model to progressively amplify reliable sources while attenuating unreliable ones.

Building on this adaptive data fusion process, the proposed MR-DRO model minimizes worst-case expected cost over a Wasserstein ambiguity set centered at a trust-weighted empirical distribution. By dynamically updating source trust and incorporating heterogeneous data in a fully nonparametric manner, the framework enhances robustness and adaptability in multi-source environments. Our  numerical experiments later demonstrate that the MR-DRO approach systematically assigns greater influence to reliable sources and achieves improved out-of-sample performance compared with existing methods.

\subsection{Summary of contributions}

The main contributions of this paper are threefold. First, we develop a DRO approach to optimize decisions under uncertainty where information about the uncertain parameters needs to be obtained from multiple data sources that have different levels of reliability. We construct the ambiguity set using nonparametric data fusion that incorporates the concept of ``trust'' and derive a  convex reformulation that preserves computational tractability. Second, we use a trust update process to simulate trust dynamics as additional data become available, using past errors to update trust. We introduce the concept of probability dominance to rationalize the existence of information sources with dominant trust. Third, through extensive computational experiments and numerical testing of instances of resource allocation and portfolio optimization, we demonstrate the effectiveness of our method and validate the derived theoretical results. 

\subsection{Structure of the paper}

The remainder of the paper is organized as follows. Section \ref{sec:model} presents the nonparametric data fusion framework and the corresponding MR-DRO model. Section \ref{sec:reformulation} develops tractable reformulations of the proposed model. Section \ref{sec:trust} introduces the dynamic trust update mechanism, characterizes the evolution of trust across data sources over time, and establishes properties related to probability dominance. Section \ref{sec:numeric} reports numerical experiments based on resource allocation and portfolio optimization problems to evaluate the computational efficiency and performance of the MR-DRO approach. Finally, Section \ref{sec:conclusion} concludes the paper and discusses directions for future research.

\paragraph{\textbf{Notation}} Throughout the paper, we use bold symbols (i.e., $\bm{u} \in \mathbb{R}^{n}$ and $\bm{\xi} \in \mathbb{R}^{m}$) to denote the vector form of a decision variable or a parameter and use bold capital letters (i.e., $\bm{A} \in \mathbb{R}^{m \times n}$) to denote a matrix. We denote the set $\{1,\ldots,N\}$ as $[N]$. The inner product of two vectors $\bm{u},\bm{v} \in \mathbb{R}^{n}$ is expressed as $\langle \bm{u},\bm{v}\rangle:= \bm{u}^\mathsf{T}\bm{v}$.

\section{Models and Solution Approaches}\label{sec:model}

We introduce the data fusion approach for handling multi-sourced reference information in Section \ref{sec:2.1}, and present the MR-DRO model in Section \ref{sec:2.2}.

\subsection{Nonparametric data fusion for multiple sources}\label{sec:2.1}
Let $M, H, I$ denote the dimension of uncertain parameters, the number of data sources, and the number of past events (up to the decision point), respectively. The uncertain parameters of our problem are captured by a random vector $\bm{\xi} \in \mathbb{R}^{M}$. Although the true distribution of $\bm{\xi}$ is inherently unknown, we rely on observations from past events to infer statistical information of $\bm{\xi}$. We use $i$ to index the events as the discrete time series of data we receive from $H$ sources, where $i \in [I]$ corresponds to past events and $i=I+1$ represents the current event at which a decision needs to be made. At each past event $i \in [I]$, we have a prediction of the uncertainty based on each information source $h \in [H]$, denoted as $\tilde{\bm{\xi}}^{(i)}_{h} \in \mathbb{R}^{M}$. Following the realization of event $i$, we observe the true value $\bm{\xi}^{(i)}_{\textrm{true}} \in \mathbb{R}^{M}$, allowing us to calculate the prediction error of source $h$, denoted as $\Delta \bm{\xi}^{(i)}_{h} = \tilde{\bm{\xi}}^{(i)}_{h} - \bm{\xi}^{(i)}_{\textrm{true}}$, for each source $h \in [H]$. 

We assume that the predictions from different sources are independent. We further assume that the prediction accuracy of each source remains consistent throughout all events, with errors from each source $h$ following an unknown distribution for all $h \in [H]$. However, it is not necessary to assume that the error distributions are identical across different sources. Therefore, utilizing the historical errors specific to each source, we can revise the predictions for the current event, and these revised predictions can serve as more reliable information references. At the beginning of the current event $I+1$, we receive predictions $\tilde{\bm{\xi}}^{(I+1)}_{h} \in \mathbb{R}^{M}$ from source $h$. Using the errors from the past $I$ events, we compute revised predictions $\hat{\bm{\xi}}_{h}$ from source $h$ for the current event as $\hat{\bm{\xi}}_{h} = [\hat{\bm{\xi}}^{(1)}_{h},\ldots,\hat{\bm{\xi}}^{(I)}_{h}]^\mathsf{T}$, $\hat{\bm{\xi}}_{h} \in \mathbb{R}^{I \times M}$. Each element $\hat{\bm{\xi}}^{(i)}_{h} \in \mathbb{R}^{M}$ for all $i \in [I]$ is calculated as $\hat{\bm{\xi}}^{(i)}_{h}=\tilde{\bm{\xi}}^{(I+1)}_{h}-\Delta \bm{\xi}^{(i)}_{h}$.

In the nonparametric data fusion, we consider trust on each data source as a weight and update the trust  after each event $i \in [I+1]$. Specifically, we assume an initial trust $\bm{t}^{(0)} = [t^{(0)}_{1},\ldots,t^{(0)}_{H}]^\mathsf{T}, \bm{t}^{(0)} \in \mathbb{R}^{H}$ before observing any data. The trust is updated to $\bm{t}^{(i)}$ after each event $i$ and applied in the subsequent event $i+1$, for all $i \in [I]$. We ensure that $\sum_{h=1}^{H}t^{(i)}_{h} = 1$ for all $i \in [I]$ and $0 \leq t^{(i)}_{h} \leq 1$ for all $i \in [I], h \in [H]$. 
Consequently, we compute the distribution $\hat{\mathbb{P}}_{HI}$ as a discrete distribution:
\begin{equation}\label{eq:empirical-distribution}
    \hat{\mathbb{P}}_{HI} := \sum_{h=1}^{H}\frac{t^{(I)}_{h}}{I}\sum_{i=1}^{I}\delta_{\hat{\bm{\xi}}^{(i)}_{h}}
\end{equation}
where $\hat{\bm{\xi}}^{(i)}_{h}$ corresponds to the $i$-th revised prediction from source $h$, and $\delta(\bm{\xi} - \hat{\bm{\xi}}^{(i)}_{h})$ is the Dirac delta function \citep{dirac1981principles}, which equals to $1$ when $\bm{\xi}$ matches the prediction $\hat{\bm{\xi}}^{(i)}_{h}$ and $0$ otherwise.

\begin{example}[Nonparametric data fusion]\label{exp:non-parametric-data-fusion}
    Consider a small example of a resource allocation problem in one region involving $H=2$ sources with only uncertain demand ($M=1$). Each source has provided predicted demands for $I=2$ past events such as $\tilde{\bm{\xi}}_{1} = [11,14]^\mathsf{T}$ and $\tilde{\bm{\xi}}_{2} = [8,14]^\mathsf{T}$. The realization of past events is given by $\bm{\xi}_{\text{true}} = [10,13]^\mathsf{T}$. With these realizations, we calculate the prediction errors for each source as $\Delta \bm{\xi}_{1}=[1,1]^\mathsf{T}$ and $\Delta \bm{\xi}_{2}=[-2,1] ^\mathsf{T}$. For the current event, we update the predicted demands for each source as $\tilde{\xi}^{(3)}_{1}=6$ and $\tilde{\xi}^{(3)}_{2}=9$, respectively. Then, we compute the revised predictions from each source for the current event’s reference as $\hat{\bm{\xi}_{1}} = [5,5]^\mathsf{T}$ and $\hat{\bm{\xi}_{2}} = [11,8]^\mathsf{T}$. Suppose that we use trust $t^{(2)}_{1} = 0.6$ in source 1 and $t^{(2)}_{2} = 0.4$ in source 2. The empirical distribution $\hat{\mathbb{P}}_{HI}$ is then composed of $H \times I = 4$ data points, with each point from the source $1$ having probability $0.3$, and each point from source $2$ having probability $0.2$.
\end{example}

\subsection{Multi-reference distributionally robust optimization (MR-DRO) model}\label{sec:2.2}
Consider decision $\bm{x} \in \mathbb{R}^{K}$ that resides in a feasible set $\mathbb{X} \subseteq \mathbb{R}^{K}$. A stochastic optimization model: 
\begin{equation}\label{eq:SP}
    \inf_{x \in \mathbb{X}}\mathbb{E}^{\mathbb{P}}\left[\ell(\bm{x},\bm{\xi})\right]
\end{equation}
seeks an optimal solution of $\bm{x}$ that minimizes the expected value of a cost-related objective function $\ell(\bm{x},\bm{\xi})$, given by $\mathbb{E}^{\mathbb{P}}\left[\ell(\bm{x},\bm{\xi})\right] = \int_{\Xi}\ell(\bm{x},\bm{\xi})\mathbb{P}(d\bm{\xi})$. Here, $\Xi \subseteq \mathbb{R}^{M}$ represents the uncertainty set, and $\ell: \mathbb{R}^{K} \times \mathbb{R}^{M} \rightarrow \mathbb{R}$ is the objective function dependent on both the decision vector $\bm{x}$ and the uncertain parameters $\bm{\xi}$.

By assumption, $\mathbb{X}$ is a closed set and the objective function $\ell(\bm{x},\bm{\xi})$ is lower semi-continuous in $\bm{x}$ for every fixed $\bm{\xi} \in \Xi$, and upper semi-continuous in $\bm{\xi}$ for every fixed $\bm{x} \in \mathbb{X}$. However, the unknown nature of the true distribution $\mathbb{P}$ prevents us from solving \eqref{eq:SP} directly. To address this issue, we formulate a multi-reference distributionally robust optimization (MR-DRO) problem:
\begin{equation}\label{eq:DRO}
        \inf_{\bm{x} \in \mathbb{X}}\left\{\sup_{\mathbb{Q} \in \mathbb{B}_{\epsilon}(\hat{\mathbb{P}}_{HI})}\mathbb{E}^{\mathbb{Q}}\left[\ell(\bm{x},\bm{\xi})\right]\right\}.
\end{equation}
Here, the ambiguity set $\mathbb{B}_{\epsilon}(\hat{\mathbb{P}}_{HI})$ is defined as:
\begin{equation}\label{eq:ambiguity-set}
   \mathbb{B}_{\epsilon}(\hat{\mathbb{P}}_{HI}) = \left\{\mathbb{Q} \in \mathcal{M}(\Xi): d_{W}(\hat{\mathbb{P}}_{HI}, \mathbb{Q}) \leq \epsilon\right\}. 
\end{equation}
which includes all distributions $\mathbb{Q}$ on the space $\mathcal{M}(\Xi)$ such that $\mathbb{E}^{\mathbb{Q}}\left[\left\|\bm{\xi}\right\|\right] = \int_{\Xi}\left\|\xi\right\|\mathbb{Q}(d\bm{\xi}) < \infty$ and their Wasserstein distance to $\hat{\mathbb{P}}_{HI}$ (obtained via \eqref{eq:empirical-distribution}) is at most $\epsilon$ ($\epsilon \geq 0$). It can be viewed as a Wasserstein ball of radius $\epsilon$ centered on the empirical distribution $\hat{\mathbb{P}}_{HI}$. The Wasserstein distance is further defined as:
\begin{equation}\label{eq:Wasserstein-metric}
    d_{W}(\hat{\mathbb{P}}_{HI}, \mathbb{Q}) = \inf \left\{ \int_{\Xi^{2}}\|\bm{\xi}-\bm{\xi}^{'}\|\Pi(d\bm{\xi},d\bm{\xi}^{'})\right\},
\end{equation}
where $\|\cdot\|$ represents an arbitrary norm on $\mathbb{R}^{M}$, and $\Pi$ is a joint distribution of $\bm{\xi}$ and $\bm{\xi}^{'}$ with marginal distributions $\mathbb{Q}$ and $\hat{\mathbb{P}}_{HI}$, respectively. 

The objective of this MR-DRO model~\eqref{eq:DRO} is to minimize the worst-case expectation of $\ell(\bm{x},\bm{\xi})$ in all possible distributions within the ambiguity set $\mathbb{B}_{\epsilon}(\hat{\mathbb{P}}_{HI})$. 
The inner maximization over probability distributions in Model~\eqref{eq:DRO} involves an optimization problem in an infinite-dimensional space, which may appear intractable. In Section \ref{sec:reformulation}, we demonstrate that under reasonable assumptions for the uncertainty set $\Xi$ and the objective function $\ell(\bm{x},\bm{\xi})$, Model~\eqref{eq:DRO} can be reformulated into computationally tractable forms.

\section{Tractable Reformulations of the MR-DRO Model}\label{sec:reformulation}
We show that the inner worst-case expectation problem in \eqref{eq:DRO} over the Wasserstein ambiguity set \eqref{eq:ambiguity-set}, expressed as:
\begin{equation}\label{eq:worst-case-expectation}
   \sup_{\mathbb{Q} \in \mathbb{B}_{\epsilon}(\hat{\mathbb{P}}_{HI})}\mathbb{E}^{\mathbb{Q}}\left[\ell(\bm{x},\bm{\xi})\right], 
\end{equation}
can be reformulated as a finite convex program for many specific forms of the  objective function $\ell(\bm{x},\bm{\xi})$. Without loss of generality, we assume that the objective function is defined as the point-wise maximum of $J$ elementary measurable functions $\ell_{j}: \mathbb{R}^{M} \times \mathbb{R}^{K} \rightarrow \mathbb{R}$, $j \in [J]$, and therefore, $\ell(\bm{x},\bm{\xi}):= \max_{j \in [J]} \left\{\ell_{j}(\bm{x},\bm{\xi})\right\}$. The following assumption on $\Xi$ and $\ell_{j}$ ensures tractability of the MR-DRO model.
\begin{assumption}[Convexity, \citet{mohajerin2018data}, Assumption 4.1]\label{asp:convexity}
    The uncertainty set $\Xi \subseteq \mathbb{R}^{M}$ is convex and closed. Each $-\ell_{j}$ is proper \citep{aliprantis2006dimensional}, convex, and lower semi-continuous in $\bm{\xi}$ for every fixed $\bm{x} \in \mathbb{X}$ and for all $j \in [J]$. In addition, $\ell_{j}$ is not identically $-\infty$ in $\Xi$ for every fixed $\bm{x} \in \mathbb{X}$ and for all $j \in [J]$.
\end{assumption}

We address the reduction of the inner worst-case expectation to a finite convex program in Section~\ref{sec:3.1}, and the reformulations for specific objective functions in Section~\ref{sec:3.2}.
\subsection{Convex reduction}\label{sec:3.1}
\begin{theorem}[Convex reduction]\label{thm:convex-reduction}
Given Assumption \ref{asp:convexity}, then for any $\epsilon \geq 0$ the inner worst-case expectation \eqref{eq:worst-case-expectation} equals the optimal objective value of the following convex program \eqref{eq:DRO-finite-convex-program}:

\begin{subequations}\label{eq:DRO-finite-convex-program}
\begin{align}
        \inf_{\lambda, s_{hi}, \bm{z}_{hij}, \bm{\nu}_{hij}}\quad &  \lambda\epsilon + \sum^{H}_{h=1}\sum^{I}_{i=1}\frac{t^{(I)}_{h}}{I} \cdot s_{hi}\label{eq:DRO-finite-convex-program-obj}\\
        \textrm{s.t.}\quad\quad\quad & [-\ell_{j}]^{*}(\bm{z}_{hij}-\bm{\nu}_{hij}) + \sigma_{\Xi}(\bm{\nu}_{hij}) - \langle \bm{z}_{hij},\hat{\bm{\xi}}^{(i)}_{h}\rangle \leq s_{hi},\label{eq:DRO-finite-convex-program-1}, \quad \forall h \in [H], i \in [I], j \in [J],\\
         & \|\bm{z}_{hij}\|_{*} \leq \lambda, \quad \forall h \in [H], i \in [I], j \in [J].\label{eq:DRO-finite-convex-program-2}
\end{align}
\end{subequations}
Here, the conjugate of a function $f$ is defined as $f^{*}(\bm{v}):= \sup_{\bm{u} \in \text{dom }f}\langle \bm{v},\bm{u}\rangle - f(\bm{u})$; notation $\|\cdot\|_{*}$ is the dual norm of the norm used in the Wasserstein distance \eqref{eq:Wasserstein-metric}, which is defined through $\|\bm{v}\|_{*} := \sup_{\|\bm{u}\| \leq 1} \langle \bm{v},\bm{u}\rangle$. The support function of $\Xi$ is defined as $\sigma_{\Xi}(\bm{u}):= \sup_{\bm{\xi} \in \Xi} \langle \bm{u}, \bm{\xi}\rangle$.
\end{theorem}
A detailed proof of Theorem \ref{thm:convex-reduction} is provided in Supplement Material Section \ref{sec:proof-1}. This result forms the foundation for further reformulations under specific structural assumptions on the loss function.

\subsection{Reformulations for certain loss functions}\label{sec:3.2}
Building on Theorem \ref{thm:convex-reduction}, we now demonstrate how the convex optimization program \eqref{eq:DRO-finite-convex-program} can be reduced to ensure computational tractability for certain loss functions.

\begin{assumption}[Piecewise affine objective functions]\label{asp:piecewise-affine-objective}
    Suppose that the uncertainty set is a polytope as $\Xi = \left\{\bm{\xi} \in \mathbb{R}^{M}: \bm{C}\bm{\xi} \leq \bm{g}\right\}$ where $\bm{C}$ is a matrix and $\bm{g}$ is a vector of appropriate dimension. In addition, for every fixed $\bm{x} \in \mathbb{X}$, each elementary function $\ell_{j}(\bm{x},\bm{\xi}) = \langle\bm{a}_{j},\bm{\xi}\rangle+b_{j}$ in the objective function $\ell(\bm{x},\bm{\xi})$ is an affine function with respect to $\bm{\xi}$, where $\bm{a}_{j} \in \mathbb{R}^{M}$ and $b_{j} \in \mathbb{R}$ are linear in $\bm{x}$. The objective function can be expressed as $\ell(\bm{x},\bm{\xi}) = \max_{j \in [J]} \left\{\langle \bm{a}_{j}, \bm{\xi}\rangle + b_{j}\right\}$ for all $j \in [J]$.
\end{assumption}
\begin{theorem}[Reformulation with piecewise affine loss functions]\label{thm:piecewise-affine}
    Given Assumption~\ref{asp:piecewise-affine-objective} then for any $\epsilon \geq 0$, the MR-DRO model \eqref{eq:DRO} is equivalent to:
\begin{subequations}\label{eq:DRO-piecewise-affine}
\begin{align}
        \inf_{\bm{x}, \lambda, s_{hi}, \bm{\gamma}_{hij}}\quad &  \lambda\epsilon + \sum^{H}_{h=1}\sum^{I}_{i=1}\frac{t^{(I)}_{h}}{I} \cdot s_{hi}\label{eq:DRO-piecewise-affine-obj}&\\
        \textrm{s.t.}\quad\quad\quad & \bm{x} \in \mathbb{X}, \label{eq:DRO-piecewise-affine-x}&\\
        & b_{j} + \langle \bm{a}_{j}, \hat{\bm{\xi}}^{(i)}_{h}\rangle + \langle \bm{\gamma}_{hij}, \bm{g} - \bm{C}\hat{\bm{\xi}}^{(i)}_{h} \rangle \leq s_{hi}, &\forall h \in [H], i \in [I], j \in [J],\label{eq:DRO-piecewise-affine-constr1}\\
        & \|\bm{C}^\mathsf{T}\bm{\gamma}_{hij} - \bm{a}_{j}\|_{*} \leq \lambda, &\forall h \in [H], i \in [I], j \in [J],\label{eq:DRO-piecewise-affine-constr2}\\
        & \bm{\gamma}_{hij} \geq 0, &\forall h \in [H], i \in [I], j \in [J].\label{eq:DRO-piecewise-affine-constr3}
\end{align}
\end{subequations}
\end{theorem}
We present the detailed proof of Theorem \ref{thm:piecewise-affine} in Supplement Material Section \ref{sec:proof-2}. This reformulation provides a computationally efficient framework for handling piecewise affine loss functions.

Note that if we use $L_{1}$-norm in the definition of Wasserstein distance \eqref{eq:Wasserstein-metric}, the dual norm here becomes $L_{\infty}$-norm. Furthermore, if the feasible set $\mathbb{X} \subseteq \mathbb{R}^{K}$ is a set of linear constraints with regard to $\bm{x}$, we have a linear programming (LP) reformulation \eqref{eq:DRO-piecewise-affine}.

\begin{assumption}[Separable affine objective functions]\label{asp:separable-affine}
    Suppose that the uncertain parameters can be separated as $\bm{\xi} = (\bm{\xi}_{1},\ldots,\bm{\xi}_{N})$, and assume that $\bm{\xi}_{n} \in \Xi_{n}$, where $\Xi_{n} \subseteq \mathbb{R}^{M}$ is nonempty and closed for any $n \in [N]$. The objective function is additively separable with respect to the structure of $\bm{\xi}$, i.e., $\ell(\bm{x},\bm{\xi}) := \sum_{n=1}^{N}\max_{j \in [J]}\left\{\ell_{nj}(\bm{x},\bm{\xi}_{n})\right\}$.
\end{assumption}

Note that such loss functions commonly appear, for example, in multi-item newsvendor problems and production planning problems under uncertain demands. We define $\|\bm{\xi}\|_{N} = \sum_{n=1}^{N}\|\bm{\xi}_{n}\|$ to measure the overall uncertainty across $N$ and also require the $\|\cdot\|$ on $\mathbb{R}^{M}$ to be the same norm we used in the definition of the Wasserstein metric \eqref{eq:Wasserstein-metric}.

\begin{theorem}[Reformulation with separable affine objective functions]\label{thm:separable-affine}
Given Assumptions~\ref{asp:piecewise-affine-objective} and~\ref{asp:separable-affine}, we have 
    $\Xi_{n} = \left\{\bm{\xi}_{n} \in \mathbb{R}^{M}: \bm{C}_{n}\bm{\xi}_{n} \leq \bm{g}_{n}\right\}$ where $\bm{C}_{n}$ is a matrix and $\bm{g}_{n}$ is a vector of appropriate dimension for all $n \in [N]$. In addition, for every fixed $\bm{x} \in \mathbb{X}$, each elementary function $\ell_{nj}(\bm{x},\bm{\xi}_{n}) = \langle\bm{a}_{nj},\bm{\xi}\rangle+b_{nj}$ is an affine function with regard to $\bm{\xi}_{n}$, where $\bm{a}_{nj} \in \mathbb{R}^{M}$ and $b_{nj} \in \mathbb{R}$ are linear in $\bm{x}$ for all $n \in [N]$ and $j \in [J]$.
    Then for any $\epsilon \geq 0$, the MR-DRO model \eqref{eq:DRO} can be reformulated as:
\begin{subequations}\label{eq:DRO-separable-affine}
\begin{align}
        \inf_{\bm{x}, \lambda, s_{hin}, \bm{\gamma}_{hijn}}\quad &  \lambda\epsilon + \sum^{H}_{h=1}\sum^{I}_{i=1}\sum^{N}_{n=1}\frac{t^{(I)}_{h}}{I} \cdot s_{hin}\label{eq:DRO-separable-affine-obj}&\\
        \textrm{s.t.}\quad\quad\quad & \bm{x} \in \mathbb{X}, \label{eq:DRO-separable-affine-x}&\\
        & b_{nj} + \langle \bm{a}_{nj}, \hat{\bm{\xi}}^{(i)}_{h,n}\rangle + \langle \bm{\gamma}_{hijn}, \bm{g}_{n} - \bm{C}_{n}\hat{\bm{\xi}}^{(i)}_{h,n} \rangle \leq s_{hin}, &\forall h \in [H], i \in [I], j \in [J], n \in [N],\label{eq:DRO-separable-affine-constr1}\\
        & \|\bm{C}_{n}^\mathsf{T}\bm{\gamma}_{hijn} - \bm{a}_{nj}\|_{*} \leq \lambda, &\forall h \in [H], i \in [I], j \in [J], n \in [N]\label{eq:DRO-separable-affine-constr2}\\
        & \bm{\gamma}_{hijn} \geq 0, &\forall h \in [H], i \in [I], j \in [J], n \in [N].\label{eq:DRO-separable-affine-constr3}
\end{align}
\end{subequations}
\end{theorem}

The proof of Theorem \ref{thm:separable-affine}, detailed in Supplement Material Section \ref{sec:proof-3}, builds on the framework established in Theorem \ref{thm:piecewise-affine} by exploiting the separability of the loss function. This allows for computational simplifications in cases where the uncertain parameters are structured as separable components.

\section{Dynamic Trust Update Process}\label{sec:trust}

We now introduce a dynamic trust update mechanism to capture how a decision maker’s trust in multiple data sources evolves over time as new information becomes available. This mechanism serves two primary objectives. First, it enhances decision accuracy and robustness by leveraging the discrete distribution $\hat{\mathbb{P}}_{HI}$—constructed through nonparametric data fusion with trust values acting as probabilistic weights—which serves as the center of the Wasserstein ambiguity set in the MR-DRO model. By dynamically updating trust levels, the framework reduces the influence of low-quality or unreliable data sources, thereby mitigating adverse effects arising from inaccurate information. Second, the adaptive adjustment of trust values based on observed performance aligns the model with realistic decision-making behavior, where confidence in different sources evolves according to their historical reliability. Increasing trust in sources that consistently yield accurate predictions while decreasing trust in poorly performing ones mirrors natural human learning processes and improves the model’s adaptability to changing environments.

One approach to modeling trust dynamics is to use realized losses from each event as a performance metric \citep{guo-etal-cdc-2024}. Under a trial-and-error framework, one may approximate the partial derivative of realized losses with respect to the trust assigned to a particular source and update trust accordingly. A negative partial derivative indicates that increasing trust in that source reduces losses, suggesting that its weight should be increased in subsequent iterations, whereas a positive derivative implies the opposite adjustment. This gradient-based approach is particularly convenient in parametric data fusion settings for constructing ambiguity sets and remains applicable even when historical prediction error data are unavailable.

In this paper, we instead propose an error-based trust update scheme that evaluates historical misalignment between each source’s predictions and the realized uncertainty outcomes, as detailed in Section \ref{sec:4.1}. Building on this idea, we generalize two classical online learning methods—the static expert algorithm and the variable-share algorithm of \citet{herbster1998tracking}—to update trust values, as presented in Sections \ref{sec:4.2} and \ref{sec:4.3}. Finally, Section \ref{sec:4.4} introduces the concept of probability dominance and establishes conditions under which sources with dominant trust can be identified.

\subsection{Min-max error trust update}\label{sec:4.1}
\begin{figure}[!ht]	
    \centering
    \includegraphics[width=0.5\textwidth]{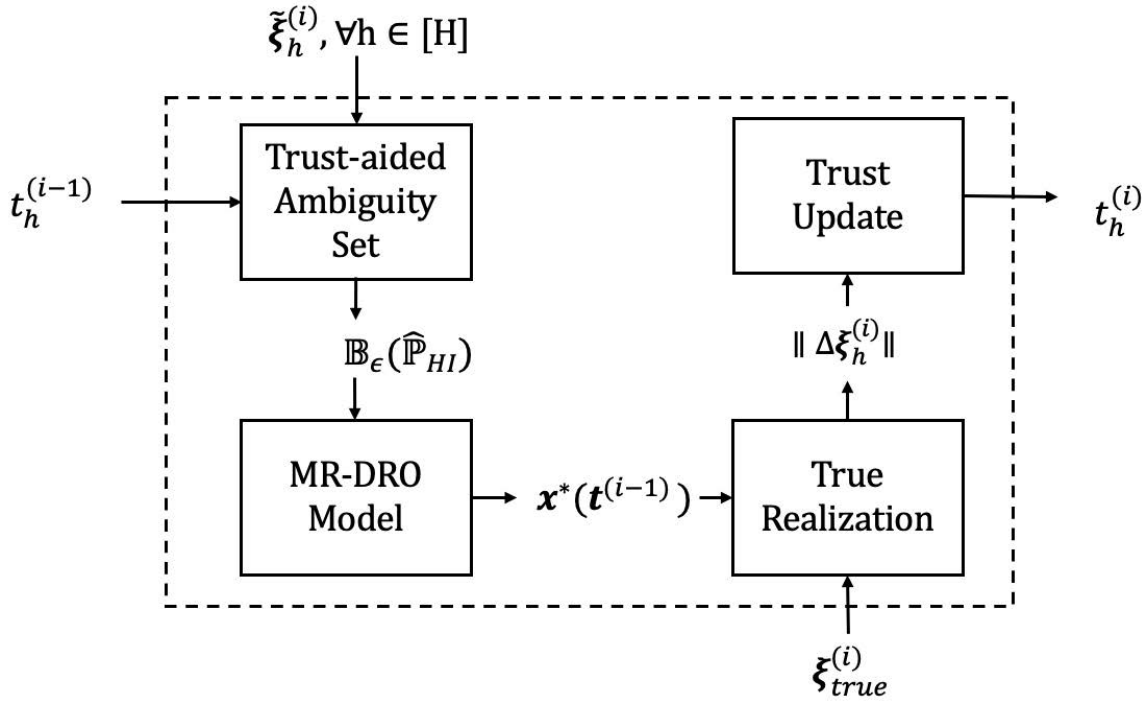}
    \caption{Illustrating trust update process with method-dependent metrics}	
    \label{fig:trust_update}
\end{figure}

We decompose each event as a four-step procedure, shown in Figure~\ref{fig:trust_update}. 
At the beginning of event $i$, we have trust $t^{(i-1)}_{h}$ for each data source $h$, for all $h \in [H]$. Then we update the ambiguity set $\mathbb{B}_{\epsilon}(\hat{\mathbb{P}}_{HI})$ and solve MR-DRO to attain an optimal solution $\bm{x}^{*}(\bm{t}^{(i-1)})$. After event $i$ ends, we obtain realization  $\bm{\xi}^{(i)}_{\textrm{true}}$ of the uncertain parameter $\bm{\xi}$. As events occur independently, each $\bm{\xi}^{(i)}_{\textrm{true}}$ is independent of the others for all $i \in [I]$.

The complete min-max error trust update process is illustrated as Algorithm~\ref{alg:min-max-error-trust-update}. Specifically, the update rule is defined as:
\begin{equation}\label{eq:min-max-error-trust-update}
    t^{(i)}_{h}=
    \begin{cases}t^{(i-1)}_{h} + \Delta t, & \textrm{if $h = \textrm{argmin}_{h \in [H]} \|\Delta \bm{\xi}^{(i)}_{h}\|$,}\\
    t^{(i-1)}_{h} - \Delta t, & \textrm{if $h = \textrm{argmax}_{h \in [H]} \|\Delta \bm{\xi}^{(i)}_{h}\|$,}\\
    t^{(i-1)}_{h}, \quad & \textrm{otherwise.}
    \end{cases}
\end{equation}
In this method, the trust update is driven by the prediction error from each source $\Delta \bm{\xi}^{(i)}_{h}$ in the most recent event, evaluated from a norm $\|\cdot\|$ on $\mathbb{R}^{M}$. For each update, we only focus on the extreme cases: a source with the minimum norm of prediction error will experience an increase in its trust value, while we will decrease the trust in the source with the maximum norm. To ensure stability, we initialize $t^{(0)}_{h} = Q/R$ for all $h \in [H]$ where $Q, R$ are positive integers with the highest common divisor $1$. The step size for the update is defined as $\Delta t = 1/qR$, where $q \in \mathbb{N}_{+}$.

\begin{algorithm}[!ht]
\caption{Min-max error trust update}\label{alg:min-max-error-trust-update}
    \DontPrintSemicolon
    \SetNoFillComment
    \KwInput{Original trust $\bm{t}^{(0)} = (t^{(0)}_{1}, \ldots, t^{(0)}_{H})$, step size $\Delta t$.}
    \For{$i=1,2,\ldots, I+1$}{
        Generate the ambiguity set with trust-weighted empirical distribution $\mathbb{B}_{\epsilon}(\hat{\mathbb{P}}_{HI})$ with $\bm{t}^{(i-1)}$.\;
        Solve MR-DRO model with $\mathbb{B}_{\epsilon}(\hat{\mathbb{P}}_{HI})$ and get the solution $\bm{x}^{*}(\bm{t}^{(i-1)})$.\;
        \For{$h=1,2,\ldots,H$}{
        Compute the last observed norm of prediction error $\|\Delta \bm{\xi}^{(i)}_{h}\|$.\;
        }
        Find $h^{+} = \textrm{argmin}_{h \in [H]} \|\Delta \bm{\xi}^{(i)}_{h}\|$.\;
        Let $(t^{\textrm{new}})^{(i-1)}_{h^{+}} = t^{(i-1)}_{h^{+}} + \Delta t$.\;
        Find $h^{-} = \textrm{argmax}_{h \in [H]} \|\Delta \bm{\xi}^{(i)}_{h}\|$.\;
        Let $(t^{\textrm{new}})^{(i-1)}_{h^{-}} = t^{(i-1)}_{h^{-}} - \Delta t$.\;
        Set $\bm{t}^{(i)} \leftarrow (\bm{t}^{\textrm{new}})^{(i-1)}$.\;
    }
\end{algorithm}

\begin{example}[Min-max error trust update]\label{exp:min-max-error-trust-update}
    Revisiting Example \ref{exp:non-parametric-data-fusion} discussed in Section~\ref{sec:2.1}, we now illustrate how we apply the min-max error trust update algorithm within the same context. Using the empirical distribution $\hat{\mathbb{P}}_{HI}$, where $H=2$ and $I=2$, we construct the ambiguity set and solve the MR-DRO model to obtain an optimal solution $x^{*}(\bm{t}^{(2)})$. We would like to update the trust to $\bm{t}^{(3)}$ for the upcoming event. Suppose that, after the current event, we obtain the realization as $\xi^{(3)}_{\text{true}} = 8$. Therefore, we obtain the prediction errors for the current event from each source as $\Delta \xi^{(3)}_{1} = -2$ and $\Delta \xi^{(3)}_{2} = 1$. Suppose that we use $L_{1}$-norm as our evaluation criterion. Since $|\Delta \xi^{(3)}_{1}| > |\Delta \xi^{(3)}_{2}|$, this indicates that the prediction of source $1$'s was less accurate than the prediction of source $2$'s for the current event. As a result, we decrease $t^{(3)}_{1}$ by a defined step size $\Delta t$ and correspondingly increase $t^{(3)}_{2}$ by $\Delta t$ for the next event.
\end{example}

\subsection{Exponential error trust update}\label{sec:4.2}

The exponential error trust update method provides an alternative framework to the min-max error trust update to dynamically adjust trust values based on observed prediction errors. Unlike the min-max approach, which focuses only on the extreme cases (the sources with the minimum and maximum prediction errors), the exponential method evaluates the prediction errors of all sources simultaneously, applying a continuous and multiplicative update mechanism.

The complete exponential error trust update process is described in Algorithm~\ref{alg:exponential-error-trust-update}. The trust update rule is defined as:
\begin{equation}\label{eq:exponential-error-trust-update}
    t^{(i)}_{h}=t^{(i-1)}_{h} \cdot e^{-\eta \cdot \|\Delta \bm{\xi}^{(i)}_{h}\|}, \quad \forall h \in [H]
\end{equation}
where $\eta > 0$ is the update rate, and $\|\cdot\|$ is a norm on $\mathbb{R}^{M}$. After computing the updated trust values, they are normalized to ensure that the sum of all trust values remains constant, preserving the interpretability of trust as weights in $\hat{\mathbb{P}}_{HI}$.

\begin{algorithm}[!ht]
\caption{Exponential error trust update}\label{alg:exponential-error-trust-update}
\DontPrintSemicolon
\KwInput{Original trust $\bm{t}^{(0)} = (t^{(0)}_{1}, \ldots, t^{(0)}_{H})$, update rate $\eta > 0$}
\For{$i=1,2,\ldots, I+1$}{
    Generate the ambiguity set with trust-weighted empirical distribution $\mathbb{B}_{\epsilon}(\hat{\mathbb{P}}_{HI})$ with $\bm{t}^{(i-1)}$.\;
    Solve MR-DRO model with $\mathbb{B}_{\epsilon}(\hat{\mathbb{P}}_{HI})$ and get the solution $\bm{x}^{*}(\bm{t}^{(i-1)})$.\;
    \For{$h=1,2,\ldots,H$}{
        Compute the last observed norm of prediction error $\|\Delta \bm{\xi}^{(i)}_{h}\|$.\;
        $(t^{\text{new}})^{(i)}_{h} = t^{(i)}_{h}\cdot e^{-\eta \cdot \|\Delta \bm{\xi}^{(i)}_{h}\|}$\;
    } 
    Normalize updated trust values $\bm{t}^{\text{new}}$\;
}
Set $\bm{t} \leftarrow \bm{t}^{\text{new}}$\;
\end{algorithm}

\begin{example}[Exponential error trust update]\label{exp:exponential-error-trust-update}
    Continuing with Example \ref{exp:min-max-error-trust-update}, we apply the exponential error trust update within the same context. Recall that the prediction errors for the current event from each source are $\Delta \xi^{(3)}_{1} = -2$ and $\Delta \xi^{(3)}_{2} = 1$. Suppose that we use $L_{1}$-norm as our evaluation criterion. With an update rate $\eta=0.5$, the trust values are updated as follows:
    \begin{equation*}
        t_{1}^{(3)} = t_{1}^{(2)} \cdot e^{-0.5\times 2}, \quad t_{2}^{(3)} = t_{2}^{(2)} \cdot e^{-0.5\times 1}.
    \end{equation*}
    After normalization, the updated trust values reflect a proportional adjustment based on the relative prediction errors of the two sources.
\end{example}

\subsection{Variable-share error trust update}\label{sec:4.3}

The variable-share error trust update method extends the exponential error trust update by introducing a redistribution mechanism that dynamically adjusts trust values based on observed prediction errors and redistributes a fraction of trust among the sources \citep{herbster1998tracking}. This approach balances penalizing poorly performing sources and rewarding better performing ones by incorporating a flexible share allocation component.

The variable share error trust update process, as described in Algorithm~\ref{alg:variable-share-error-trust-update}, begins with a multiplicative trust adjustment to obtain the intermediate trust values $(t^{\prime})^{(i)}_{h}$ for all $h \in [H]$ based on the observed prediction error for each source, similar to the exponential method. Next, we introduce a pooling mechanism, where a portion of the adjusted trust $W_{\text{trust}}$ is redistributed. This redistribution is governed by a pre-selected parameter value $\beta \in (0,1]$, which determines the fraction of the trust pool available for reallocation. The redistribution formula ensures that trust values are influenced by both individual performance and contributions from other sources:
\begin{align}\label{eq:variable-share-error-trust-update}
    & W^{(i)}_{\text{trust}} := \sum_{h=1}^{H} (1 - (1-\beta)^{\|\Delta \bm{\xi}^{(i)}_{h}\|}) (t')^{(i)}_{h},\nonumber\\
    & (t^{\text{new}})^{(i)}_{h} = (1 - \beta)^{\|\Delta \bm{\xi}^{(i)}_{h}\|} \cdot (t^{\prime})^{(i)}_{h} + \frac{1}{H-1} \left(W^{(i)}_{\text{trust}} - (1 - (1 - \beta)^{\|\Delta \bm{\xi}^{(i)}_{h}\|}) \cdot (t^{\prime})^{(i)}_{h} \right).
\end{align}
This redistribution mechanism improves collaboration among sources and ensures that trust values evolve more inclusively. Since before a source is the most accurate starting after a certain number of events, the exponential term associated with its prediction deviation may be arbitrarily large, and thus its trust may become arbitrarily small. The variable-share error trust update algorithm ensures that small trust values can be recovered quickly by the redistribution mechanism, making it more effective when the prediction precisions of different sources change during the trust update process.  

\begin{algorithm}[!ht]
\caption{Variable-share error trust update}\label{alg:variable-share-error-trust-update}
\DontPrintSemicolon
\KwInput{Original trust $\bm{t}^{(0)} = (t^{(0)}_{1}, \ldots, t^{(0)}_{H})$, share portion $\beta \in (0, 1]$, update rate $\eta > 0$}
\For{$i=1,2,\ldots, I+1$}{
    Generate the ambiguity set with trust-weighted empirical distribution $\mathbb{B}_{\epsilon}(\hat{\mathbb{P}}_{HI})$ with $\bm{t}^{(i-1)}$.\;
    Solve MR-DRO model with $\mathbb{B}_{\epsilon}(\hat{\mathbb{P}}_{HI})$ and get the solution $\bm{x}^{*}(\bm{t}^{(i-1)})$.\;
    \For{$h=1,2,\ldots,H$}{
        Compute the last observed norm of prediction error $\|\Delta \bm{\xi}^{(i)}_{h}\|$.\;
        Compute intermediate trust: $(t^{\prime})^{(i)}_{h} = t^{(i)}_{h}\cdot e^{-\eta \cdot \|\Delta \bm{\xi}^{(i)}_{h}\|}$\; 
    }
    Compute the pool: $W^{(i)}_{\text{trust}} := \sum_{h \in [H]} (1 - (1 - \beta)^{\bm{v}^{(i)}_{h}}) \cdot (t^{\prime})^{(i)}_{h}$\;
    \For{$h=1,2,\ldots,H$}{
        $(t^{\text{new}})^{(i)}_{h} = (1 - \beta)^{\|\Delta \bm{\xi}^{(i)}_{h}\|} \cdot (t^{\prime})^{(i)}_{h} + \frac{1}{H-1} \left(W^{(i)}_{\text{trust}} - (1 - (1 - \beta)^{\|\Delta \bm{\xi}^{(i)}_{h}\|}) \cdot (t^{\prime})^{(i)}_{h} \right)$\;
    }  
    Normalize updated trust values $\bm{t}^{\text{new}}$\;
}
Set $\bm{t} \leftarrow \bm{t}^{\text{new}}$\;
\end{algorithm}

\begin{example}[Variable-share error trust update]\label{exp:variable-share-error-trust-update}
    Building on Example \ref{exp:exponential-error-trust-update}, we illustrate how to apply the update to the trust of variable-share error. Using $L_{1}$-norm for error evaluation, an update rate $\eta=0.5$ and a sharing portion $\beta=0.01$, we first compute intermediate trust values as:
    \begin{equation*}
        (t^{'})_{1}^{(3)} = t_{1}^{(2)} \cdot e^{-0.5\times 2}, \quad (t^{'})_{2}^{(3)} = t_{2}^{(2)} \cdot e^{-0.5\times 1}.
    \end{equation*}
    Then, we construct the pool $W^{(i)}_{\text{trust}}$ and redistribute the pooled trust among all sources. Finally, we normalize the trust values that will be used as weights in the construction of the ambiguity set for the next event.
\end{example}

\subsection{Analyzing trust dynamics}\label{sec:4.4}

We are particularly interested in analyzing the dynamics of trust assigned to different sources under specific trust update algorithms and understanding whether and when dominant trust exists. We focus on two methods: the min-max error trust update in Section \ref{sec:4.1} and the exponential error trust update in Section \ref{sec:4.2}. 
We exclude the variable-share error update method in Section \ref{sec:4.3} from the analysis because its primary objective is to rapidly adjust trust levels when relative accuracy relationships between sources change. Consequently, the concept of dominant trust is less relevant to the variable-share error update method.

\subsubsection{Dominant trust under min-max error trust update}
\label{sec:4.4.1}
We consider the concept of probability dominance \citep{wrather1982probability} to rationalize the assignment of dominant trust (if it exists) to one source over others. 
Let $\Omega$ denote a collection of random variables corresponding to the deviations of the predicted values from the true realizations, specifically the norm of the prediction errors. For simplicity, we first examine the case where there are only two data sources, $h_Y$ and $h_Z$, that is, $H=2$. Let $Y = \|\Delta \bm{\xi}^{(i)}_{h_{Y}}\|$ represent the norm of the prediction errors of source $h_{Y}$, and $Z = \|\Delta \bm{\xi}^{(i)}_{h_{Z}}\|$ for the source $h_{Z}$, over multiple events $i = 1,2,\ldots$. The cumulative distribution functions of $Y$ and $Z$ are well defined and are denoted by $F_{Y}(a)$ and $F_{Z}(a)$, respectively. In this framework, a source that consistently exhibits smaller deviations is considered to be more reliable, and the trust assigned to this source is updated accordingly. The concept of probability dominance formalizes this comparison.

\begin{definition}[Probability dominance, \citet{wrather1982probability}]\label{def:definition1}
    Given two random outcomes, $Y, Z \in \Omega$, we say that $Y$ dominates $Z$ with probability $\beta \geq 0.5$, denoted by $Y\beta Z$, if and only if \[\mathbb{P}[Y < Z] \geq \beta.\]
\end{definition}

Based on this definition, we argue that $Y\beta Z$ is a necessary condition for fully trusting one source over another in a pairwise comparison. This implies that a decision maker's trust in source $1$ with deviations $Y$ will converge to $1$ after a finite number of events. The detailed proof of the following theorem can be found in Supplement Material Section \ref{sec:proof-4}.

\begin{theorem}[Dominant trust under min-max error trust update for two sources]\label{thm:theorem4}
    If $Y\beta Z$, then after a finite number of events, the trust in source $h_{Y}$ with deviation $Y$ will go to $1$ and fluctuates within a small interval around this maximum level.
\end{theorem}

Furthermore, one can generalize the above result to cases with multiple data sources, i.e., $H \geq 3$, stated in Theorem \ref{thm:theorem5} below. The detailed proof is provided in Supplement Material Section \ref{sec:proof-5}.

\begin{theorem}[Dominant trust under min-max error trust update for multiple sources]\label{thm:theorem5} 
    Suppose that $Y$ dominates the deviations of all other sources $Z_{h}$ in pairwise comparisons with probability $\beta \geq 0.5$, $h \neq h_{Y}, h \in H$. Specifically, $\mathbb{P}[Y < Z_{h}] \geq \beta$. Then, under the min-max error trust update Algorithm~\ref{alg:min-max-error-trust-update}, the trust in the source with deviation $Y$ will converge to $1$ after a finite number of events and fluctuate within a small interval around this maximum level.
\end{theorem}

We further investigate the relationship between first-order stochastic dominance (FSD) and probability dominance when $Y$ and $Z$ are independent, with $F_{Y}(a)$ and $F_{Z}(a)$ continuous.

\begin{theorem}[Probability dominance and FSD, \citet{wrather1982probability}]\label{thm:theorem6} 
    $Y$ and $Z$ are independent, with well-defined $F_{Y}(a)$ and $F_{Z}(a)$. If $Y$ stochastically dominates $Z$ in the first degree, then $\mathbb{P}[Y < Z] > \beta$ for $\beta=0.5$.
\end{theorem}

We present the proof of Theorem~\ref{thm:theorem6} in Supplement Material Section~\ref{sec:proof-6}.

\subsubsection{Dominant trust under exponential error trust update}

In this section, we establish the conditions for dominant trust when using the exponential error trust update Algorithm~\ref{alg:exponential-error-trust-update}. Unlike the min-max error trust update Algorithm~\ref{alg:min-max-error-trust-update}, this approach does not require probability dominance. Instead, it relies on the consistency of smaller expected prediction errors over time. The detailed proof of the following theorem can be found in Supplement Material Section \ref{sec:proof-7}.

\begin{theorem}[Dominant trust under exponential error trust update for multiple sources]\label{thm:theorem7} 
    Suppose that among $H$ sources, source $h_{Y}$ consistently produces smaller expected prediction errors than any other source. Specifically, assume there exists a constant $\zeta >0$ such that for all $h \neq h_{Y}$,
\begin{equation*}
    \mathbb{E}\left[\|\Delta \xi^{(i)}_{h_{Y}}\|\right] + \zeta \leq \mathbb{E}\left[\|\Delta \xi^{(i)}_{h}\|\right],
\end{equation*}
    for all events $i=1,2,\ldots$. Also, assume that the errors $\|\Delta\xi^{(i)}_{h}\|$ are independent over time and bounded within the interval $[0,L_{\text{max}}]$. Then, using the exponential error trust update Algorithm~\ref{alg:exponential-error-trust-update}, the trust $t^{(i)}_{h_{Y}}$ will converge to $1$ as $i \rightarrow \infty$.
\end{theorem}

\section{Computational Results}
\label{sec:numeric}

We validate the theoretical results of the MR-DRO approach in two distinct contexts, each characterized by its unique appearance of uncertainty, i.e., a resource allocation problem and a portfolio optimization problem, respectively. In resource allocation, the problem is formulated with uncertain demand appearing in the constraints. In portfolio optimization, the uncertain parameter is embedded within the objective function, specifically in the returns of the assets.

All linear programming models (i.e., reformulations of the MR-DRO model) are optimized by Gurobi 9.5.2. The algorithm for trust update is implemented in Python 3.9.12. All numerical tests are conducted on a PC with 16 GB RAM and an Apple M1 Pro chip.

\subsection{Resource allocation with uncertain demand}\label{sec:5.1}
\subsubsection{MR-DRO for a resource allocation problem}\label{sec:5.1.1}
We optimize the allocation of resources in $K_{\text{r}}$ regions, constrained by an overall budget $B > 0$. The demands $\bm{d} = [d_{1},\ldots,d_{K_{\text{r}}}]^\mathsf{T}$, $\bm{d} \in \mathbb{R}^{K_{\text{r}}}_{+}$ in all regions are uncertain and they can be estimated using predictions from multiple $h=1,\ldots,H$ sources. Let $\bm{c}^{u} = [c^{u}_{1},\ldots,c^{u}_{K_{\text{r}}}]^\mathsf{T}$ be the unit penalty cost of unmet demand and $\bm{c}^{o} = [{c}^{o}_{1},\ldots,{c}^{o}_{K_{\text{r}}}]^\mathsf{T}$ be that of over-served demand in each region. We specify the decision vector $\bm{x}=[x_{1},\ldots,x_{K_{\text{r}}}]^\mathsf{T}$ with a feasible region $\mathbb{X} = \{\bm{x} \in \mathbb{R}^{K_{\text{r}}}_{+}: \sum_{k=1}^{K_{\text{r}}}x_{k} \leq B\}$, where $x_{k}$ indicates the amount of resource allocated to region $k$, for each $k \in [K_{\text{r}}]$. For any vector $\bm{u} \in \mathbb{R}^{M}$, we define the function $(\bm{u})^{+} : = [\max\{u_{1},0\},\ldots,\max\{u_{M},0\}]^\mathsf{T}$. 
Following the structure of MR-DRO given in Section \ref{sec:2.2}, we formulate:  
\begin{equation}
    J_{\textrm{resource}}^{*} =\inf_{\bm{x} \in \mathbb{X}} \left\{\sup_{\mathbb{Q} \in \mathbb{B}_{\epsilon}(\hat{\mathbb{P}}_{HI})}\mathbb{E}^{\mathbb{Q}}\left[(\bm{c}^{u})^{\mathsf{T}}(\bm{d}-\bm{x})^{+}+(\bm{c}^{o})^{\mathsf{T}}(\bm{x}-\bm{d})^{+}\right]\right\},
\end{equation}
which minimizes the expected losses caused by undesirable allocation. By expressing the loss function as the additively separable form, we obtain:
\begin{align}
\label{eq:MRDRO:resource-allocation}
    J_{\textrm{resource}}^{*} &= \inf_{\bm{x} \in \mathbb{X}}\left\{\sup_{\mathbb{Q} \in \mathbb{B}_{\epsilon}(\hat{\mathbb{P}}_{HI})}\mathbb{E}^{\mathbb{Q}}\left[\sum_{k=1}^{K_{\text{r}}}\max\left\{c^{u}_{k}(d_{k}-x_{k}),c^{o}_{k}(x_{k}-d_{k})\right\}\right]\right\}\nonumber\\
    &= \inf_{\bm{x} \in \mathbb{X}}\left\{\sup_{\mathbb{Q} \in \mathbb{B}_{\epsilon}(\hat{\mathbb{P}}_{HI})}\mathbb{E}^{\mathbb{Q}}\left[\sum_{k=1}^{K_{\text{r}}}\max_{j \in [J]}(a_{jk}d_{k}+b_{jk})\right]\right\},
\end{align}
where $J=2$, $a_{1k} = c^{u}_{k}$, $a_{2k} = -c^{o}_{k}$, $b_{1k} = -c^{u}_{k}x_{k}$, $b_{2k}=c^{o}_{k}x_{k}$. We further assume that the predictions from each source $h$ are independent in different regions. As a result, we have separate trust after event $i$ in each region $k$ for each source $h$, denoted as $t^{(i)}_{h,k}$, for all $h \in [H]$, $i \in [I]$, $k \in [K]$.
Then by Theorem~\ref{thm:separable-affine}, we can solve the following convex program \eqref{eq:DRO-resource} as an equivalent reformulation of \eqref{eq:MRDRO:resource-allocation} for an optimal solution of the resource allocation plan $\bm{x}$:
\begin{subequations}\label{eq:DRO-resource}
\begin{align}
    \inf_{\bm{x}, \lambda, s_{hik}, \bm{\gamma}_{hijk}}\quad & \lambda\epsilon + \sum^{H}_{h=1}\sum^{I}_{i=1}\sum^{K_{\text{r}}}_{k=1}\frac{t^{(I)}_{h,k}}{I} \cdot s_{hik}\label{eq:DRO-resource-obj} & \\
    \textrm{s.t.}\quad\quad\quad & \sum_{k=1}^{K_{\text{r}}}x_{k} \leq B, \label{eq:DRO-resource-x} & \\
    & b_{jk} + a_{jk}\hat{d}^{(i)}_{h,k} + \langle \bm{\gamma}_{hijk}, \bm{g}_{k} - \bm{C}_{k}\hat{d}^{(i)}_{h,k} \rangle \leq s_{hik}, & \forall i \in [I], j \in [J], k \in [K_{\text{r}}], h \in [H],\label{eq:DRO-resource-constr1}\\
    & \|\bm{C}_{k}^\mathsf{T}\bm{\gamma}_{hijk} - a_{jk}\|_{*} \leq \lambda, & \forall i \in [I], j \in [J], k \in [K_{\text{r}}], h \in [H],\label{eq:DRO-resource-constr2}\\
    & x_{k} \geq 0, \bm{\gamma}_{hijk} \geq 0, & \forall i \in [I], j \in [J], k \in [K_{\text{r}}], h \in [H].\label{eq:DRO-resource-constr3}
\end{align}
\end{subequations}

\subsubsection{Baseline settings for the resource allocation problem}\label{sec:5.1.2}
\begin{figure}[htbp!]
    \centering
    \includegraphics[width=0.9\textwidth]{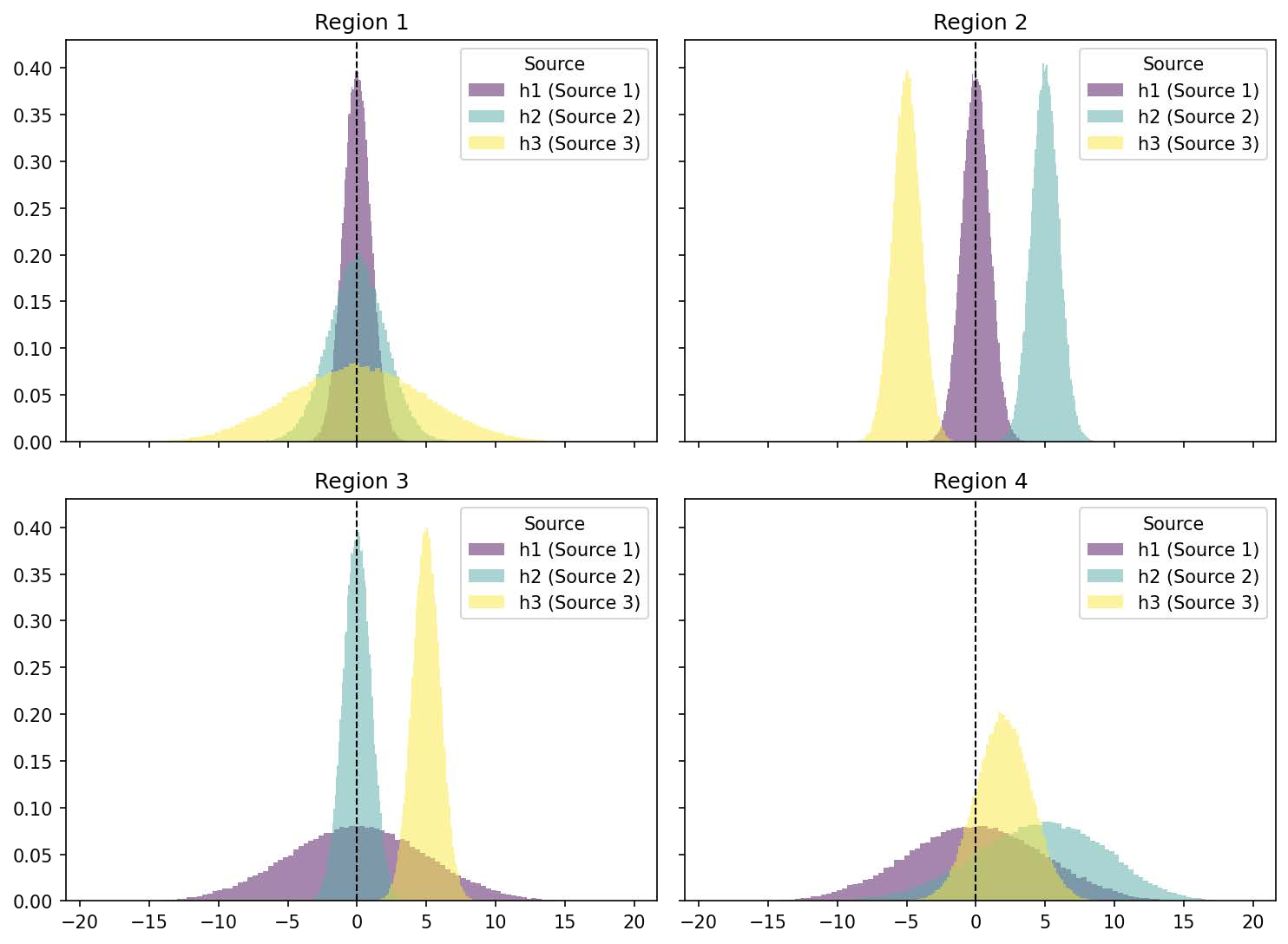}
    \caption{Prediction error distributions with the baseline setting for the resource allocation problem}
    \label{fig:error-distribution}
\end{figure}

We consider a baseline case with $K_{\text{r}}=4$, $H=3$, and $I=200$. We set $\epsilon = 0.01$ as the radius of the ambiguity set $\mathcal{P}$. The unit penalty costs for unmet and over-served demand are $\bm{c}^{u} = (\$5000,\$5000,\$5000,\$5000)^\mathsf{T}$ and $\bm{c}^{o} = (\$1000,\$1000,\$1000,\$1000)^\mathsf{T}$, respectively. The resource budget is $B = 200$. The true demand $\bm{d}^{(i)}_{\textrm{true}}$ for event $i$ is a vector with each element sampled from a uniform distribution $U(10,20)$. The prediction at region $k$ from source $h$ for event $i$ is sampled from a truncated normal distribution with mean $(d^{(i)}_{\textrm{true},k}+\mu_{h,k})$ and variance $\sigma^{2}_{h,k}$ within interval $[0,30]$. Specifically, we set $\bm{\mu} = (\bm{\mu}_{1},\bm{\mu}_{2},\bm{\mu}_{3})^\mathsf{T} = ((0,0,0,0)^\mathsf{T},(0,5,0,5)^\mathsf{T},(0,-5,5,2)^\mathsf{T})^\mathsf{T}$ and $\bm{\sigma} = (\bm{\sigma}_{1},\bm{\sigma}_{2},\bm{\sigma}_{3})^\mathsf{T} = ((1,1, 5,5)^\mathsf{T},(2,1,1,5)^\mathsf{T},(5,1,1,2)^\mathsf{T})^\mathsf{T}$. The relationship between the prediction error distributions of each source and the true realization value is demonstrated in Figure~\ref{fig:error-distribution}. The step size for min-max error trust update is $\Delta t=0.01$. The update rate is $\eta=0.5$ for both exponential error trust update and variable-share error trust update methods, and the share portion is $\beta = 0.01$ in the variable-share trust update method. We conduct $30$ trials with distinct random seeds for each trust update algorithm. 

\subsubsection{Baseline results for the resource allocation problem}\label{sec:5.1.3}

\begin{figure}[htbp]
\centering
\begin{subfigure}{\textwidth}
  \centering
  \includegraphics[width=\textwidth]{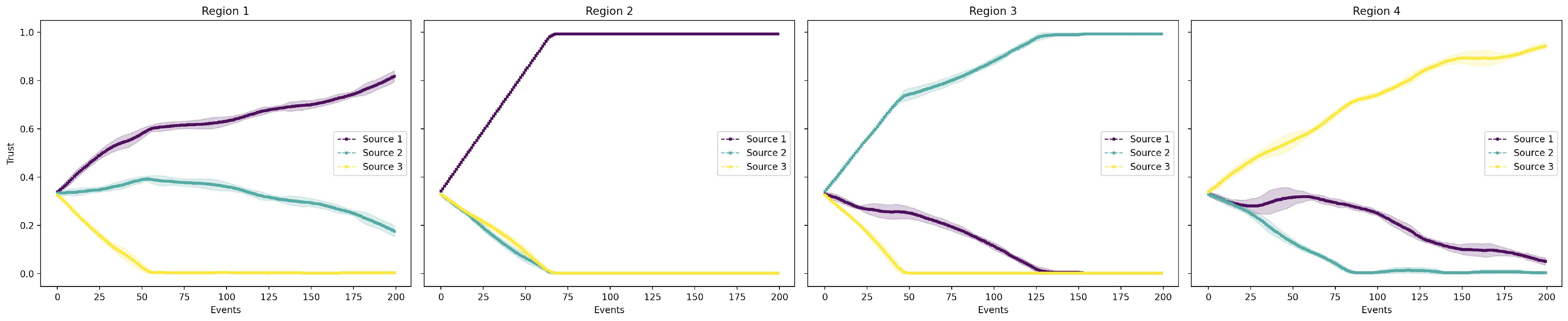}
  \caption{Min-max error trust update}
  \label{fig:resource_trust_update_min_max}
\end{subfigure}
\vspace{0.5em}  
\begin{subfigure}{\textwidth}
  \centering
  \includegraphics[width=\textwidth]{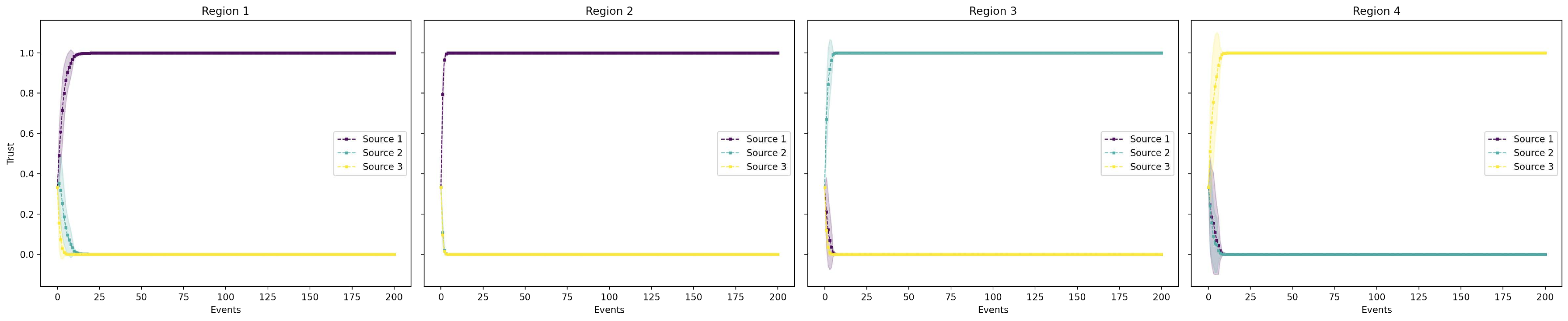}
  \caption{Exponential error trust update}
  \label{fig:resource_trust_update_exponential}
\end{subfigure}
\vspace{0.5em}  
\begin{subfigure}{\textwidth}
  \centering
  \includegraphics[width=\textwidth]{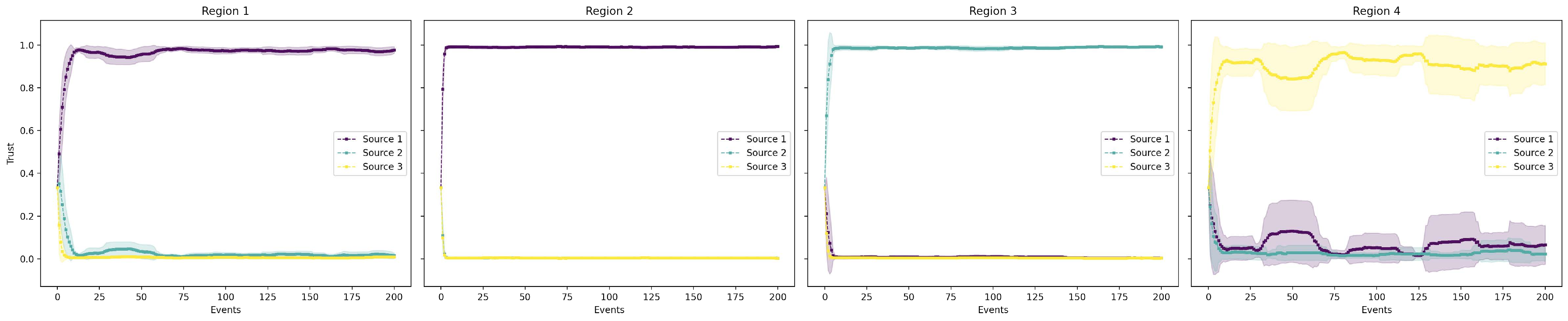}
  \caption{Variable-share error trust update}
  \label{fig:resource_trust_update_variable_share}
\end{subfigure}
\caption{Trust Update Process with 
(\protect\subref{fig:resource_trust_update_min_max}) Min-max error trust update; 
(\protect\subref{fig:resource_trust_update_exponential}) Exponential error trust update; 
(\protect\subref{fig:resource_trust_update_variable_share}) Variable-share error trust update
in the baseline setting for the resource allocation problem 
($K_{\text{r}} = 4$, $H = 3$, $I = 200$).}
\label{fig:resource_trust_update_baseline}
\end{figure}

With the baseline settings and trust update algorithms, we analyze the results after $I=200$ events, as reported in Table~\ref{tab:performance-resource} and Figure~\ref{fig:resource_trust_update_baseline}. In Table~\ref{tab:performance-resource}, the columns labeled ``Model'', ``Obj.'', ``Loss'' and ``Time'' represent the model (with the corresponding trust update algorithm) used to solve a resource allocation plan, the objective value (in thousands of dollars), the average loss (in thousands of dollars), and the total computation time (in seconds), respectively. The last three columns report the mean results across $30$ trials, with variability expressed as the mean $\pm$ standard deviation. In particular, the objective value and average loss measure different aspects of performance. The objective value reflects the worst-case expectation over a set of distributions within a Wasserstein ball centered on the trust-weighted empirical distribution (incorporating data from all past events).  In contrast, the average loss quantifies the mean actual loss observed after uncertainty is realized, based on the solution from the model. In Figure~\ref{fig:resource_trust_update_baseline}, the solid lines represent the mean trust values for each source across the $30$ trials, the shaded regions represent the range of mean $\pm$ standard deviation, illustrating the variability in trust updates between trials.
 
\begin{table}[htbp!]
\centering
\caption{Performances of different trust update methods with the baseline setting for the resource allocation problem ($K_{\text{r}} = 4$, $H = 3$ and $I = 200$)}
\begin{tabular}{c c c c} 
\toprule
Model &  Obj. (*\$1000) & Loss (*\$1000) & Time (sec.)\\
\midrule
MR-DRO (Min-max) & $8.827 \pm 1.319$ & $8.331 \pm 0.068$ & $120.907 \pm 33.557$\\
MR-DRO (Exponential) & $7.216 \pm 0.008$ & $7.390 \pm 0.095$ & $111.219 \pm 35.420$\\
MR-DRO (Variable-share) & $8.342 \pm 1.129$ & $7.332 \pm 0.097$ & $127.950 \pm 35.194$\\
DRO (h1) & $13.839 \pm 2.266$ & $15.782 \pm 0.231$ & $104.111 \pm 33.179$\\
DRO (h2) & $12.874 \pm 0.629$ & $13.789 \pm 0.123$ & $103.652 \pm 33.468$\\
DRO (h3) & $11.891 \pm 1.260$ & $13.333 \pm 0.173$ & $103.961 \pm 33.428$\\
\bottomrule
\end{tabular}
\label{tab:performance-resource}
\end{table}

The results demonstrate that the MR-DRO models exhibit the best performance among all the DRO models in both objective value and average loss. The superior performance of the MR-DRO models can be attributed to its ability to mitigate the impact of prediction errors from different sources by leveraging trust values obtained from the trust update process. This highlights the advantages of incorporating reference information from multiple sources and trust updates in the decision-making process. Furthermore, MR-DRO models with different trust update methods exhibit comparable performance in the baseline setting. 

We also perform an out-of-sample performance analysis with the $\bm{x}^{*}(\bm{t}^{(I)})$ solved by MR-DRO models and the solutions $\bm{x}^{*}_{h}$ solved by DRO models with information from source $h$, $\forall h \in [H]$. The out-of-sample test is conducted with an empirical distribution constructed by $|\Psi|$ events with each event $\psi \in \Psi$ having an equal probability $p^{\psi} = 1/|\Psi|$. The out-of-sample performance is evaluated as:
\begin{equation}
    J_{\textrm{resource}}(\bm{x}^{*}) :=\frac{1}{|\Psi|}\sum_{\psi \in \Psi}\left[(\bm{c}^{u})^{\mathsf{T}}(\bm{d}^{\psi}-\bm{x}^{*})^{+}+(\bm{c}^{o})^{\mathsf{T}}(\bm{x}^{*}-\bm{d}^{\psi})^{+}\right].\nonumber
\end{equation}

The results from the out-of-sample test, summarized in Table~\ref{tab:out-of-sample-performance-resource}, consistently highlight the best performance of MR-DRO models compared to DRO models with single-source reference information, especially when no individual source consistently provides the best predictions across all regions.

\begin{table}[htbp!]
\centering
\caption{Out-of-sample performances of different methods for the resource allocation problem ($K_{\text{r}} = 4$ and $H = 3$)}
\begin{tabular}{c c c} 
\toprule
$|\Psi|$ & Model & Loss (*\$1000)\\
\midrule
\multirow{ 5}{*}{$40$} & MR-DRO (Min-max) & $11.992 \pm 3.765$\\
 & MR-DRO (Exponential) & $12.866 \pm 4.042$\\
 & MR-DRO (Variable-share) & $12.460 \pm 4.029$\\
 & DRO (h1) & $25.675 \pm 8.188$\\
 & DRO (h2) & $44.803 \pm 16.536$\\
 & DRO (h3) & $43.208 \pm 8.235$\\
\bottomrule
\end{tabular}
\label{tab:out-of-sample-performance-resource}
\end{table}

\subsubsection{Sensitivity analysis for the resource allocation problem}\label{sec:5.1.4}
\paragraph{Varying budget $B$.}

    We maintain all other parameters as specified in the baseline setting and set the budget to $B=60$ to investigate the performance of the MR-DRO model and other methods under a limited resource allocation budget. With the same random seed setting, the trust update process remains consistent with that shown in Figure~\ref{fig:resource_trust_update_baseline}, as the prediction errors of each source align with the baseline setting. Performance results are presented in Table~\ref{tab:performance-resource-budget}, which indicates that regardless of whether the budget is always sufficient or occasionally insufficient, MR-DRO models generally outperform DRO models using single-source reference information, in terms of both objective value and average losses.

\begin{table}[htbp!]
\centering
\caption{Performances of different methods for the resource allocation problem with varying budget ($K_{\text{r}} = 4$, $H = 3$, $I = 200$ and $B = 60$)}
\begin{tabular}{c c c c c} 
\toprule
Model & Obj. (*\$1000) &  Loss (*\$1000) & Time (sec.)\\
\midrule
MR-DRO (Min-max) & $19.028 \pm 14.743$ & $15.766 \pm 0.413$ & $127.955 \pm 36.019$\\
MR-DRO (Exponential) & $17.234 \pm 15.261$ & $16.351 \pm 0.382$ & $116.484 \pm 36.533$\\
MR-DRO (Variable-share) & $18.178 \pm 14.822$ & $16.281 \pm 0.368$ & $180.354 \pm 43.576$\\
DRO (h1) & $30.778 \pm 22.675$ & $25.076 \pm 0.481$ & $110.348 \pm 35.852$\\
DRO (h2) & $26.536 \pm 19.188$ & $21.884 \pm 0.399$ & $109.940 \pm 36.121$\\
DRO (h3) & $22.620 \pm 12.010$ & $21.328 \pm 0.487$ & $110.071 \pm 36.016$\\
\bottomrule
\end{tabular}
\label{tab:performance-resource-budget}
\end{table}

\paragraph{Varying number of events $I$.}

\begin{figure}[htbp!]
    \centering
    \includegraphics[width=0.7\textwidth]{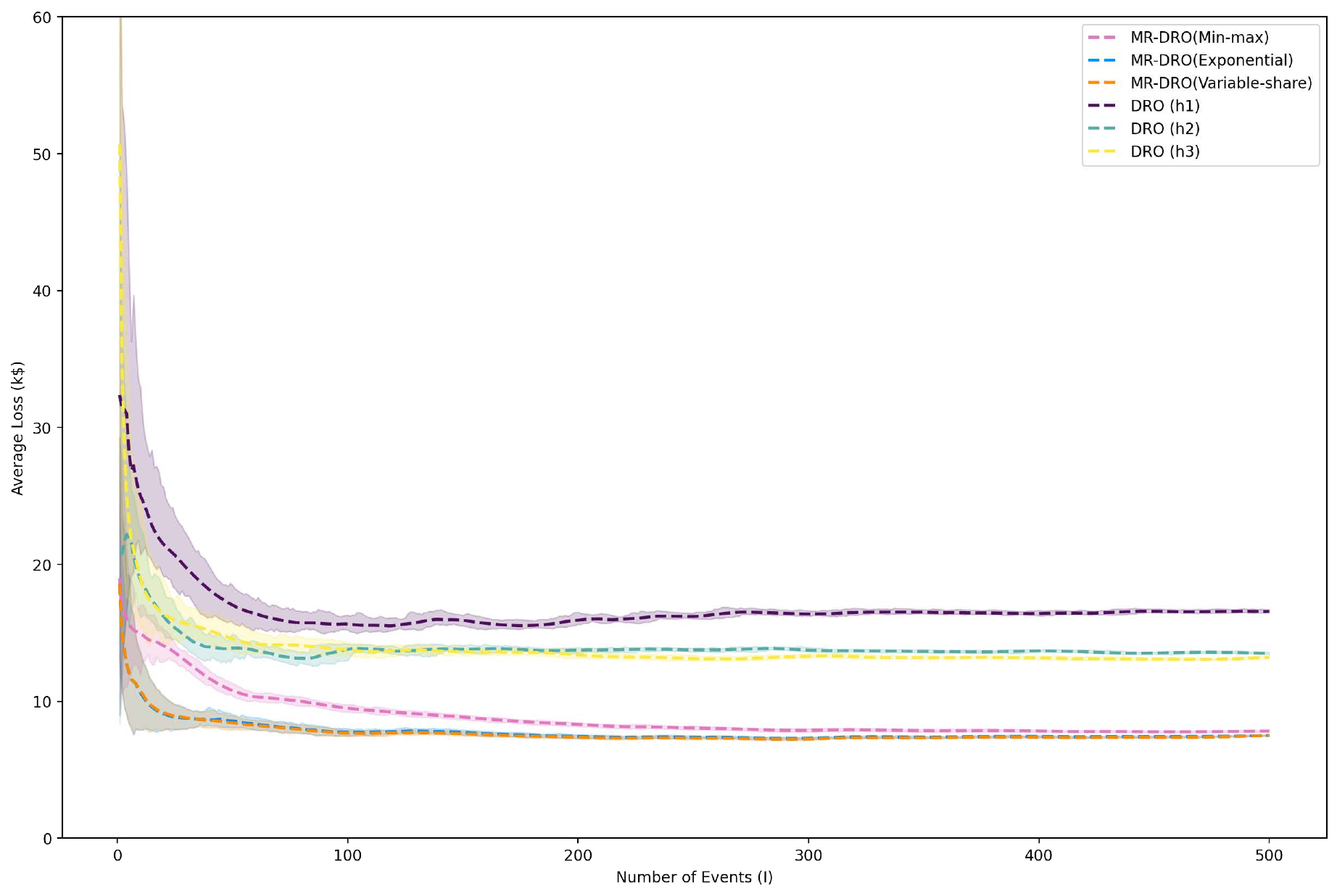}
    \caption{Average loss of different methods as number of events increases for the resource allocation problem}
    \label{fig:varying-I}
\end{figure}

    We conduct experiments to assess whether variations in $I$ affect the performance of different models. Figure~\ref{fig:varying-I} demonstrates the range of average loss across $30$ trials as the number of events increases. It shows that the average loss of the MR-DRO models initially decreases and then stabilizes, indicating that, as more events are observed, the models effectively incorporate the growing dataset and reach steady-state performance. It also suggests that a finite number of past events would be enough for us to update trust correspondingly and obtain a satisfactory solution. 

    Meanwhile, among all dynamic trust update mechanisms, the exponential and variable-share error trust update methods converge to a stable phase more quickly than the min-max error trust update method, requiring fewer iterations to reach steady performance. This suggests that these methods may be more efficient in adapting to uncertainty when the number of events is limited or when computational efficiency is a priority. In contrast, the min-max error trust update method, while still achieving strong performance, exhibits slower convergence, with the average loss stabilizing after a larger number of events.

\paragraph{Varying distribution types.}

\begin{figure}[htbp!]
    \centering
    \includegraphics[width=0.9\textwidth]{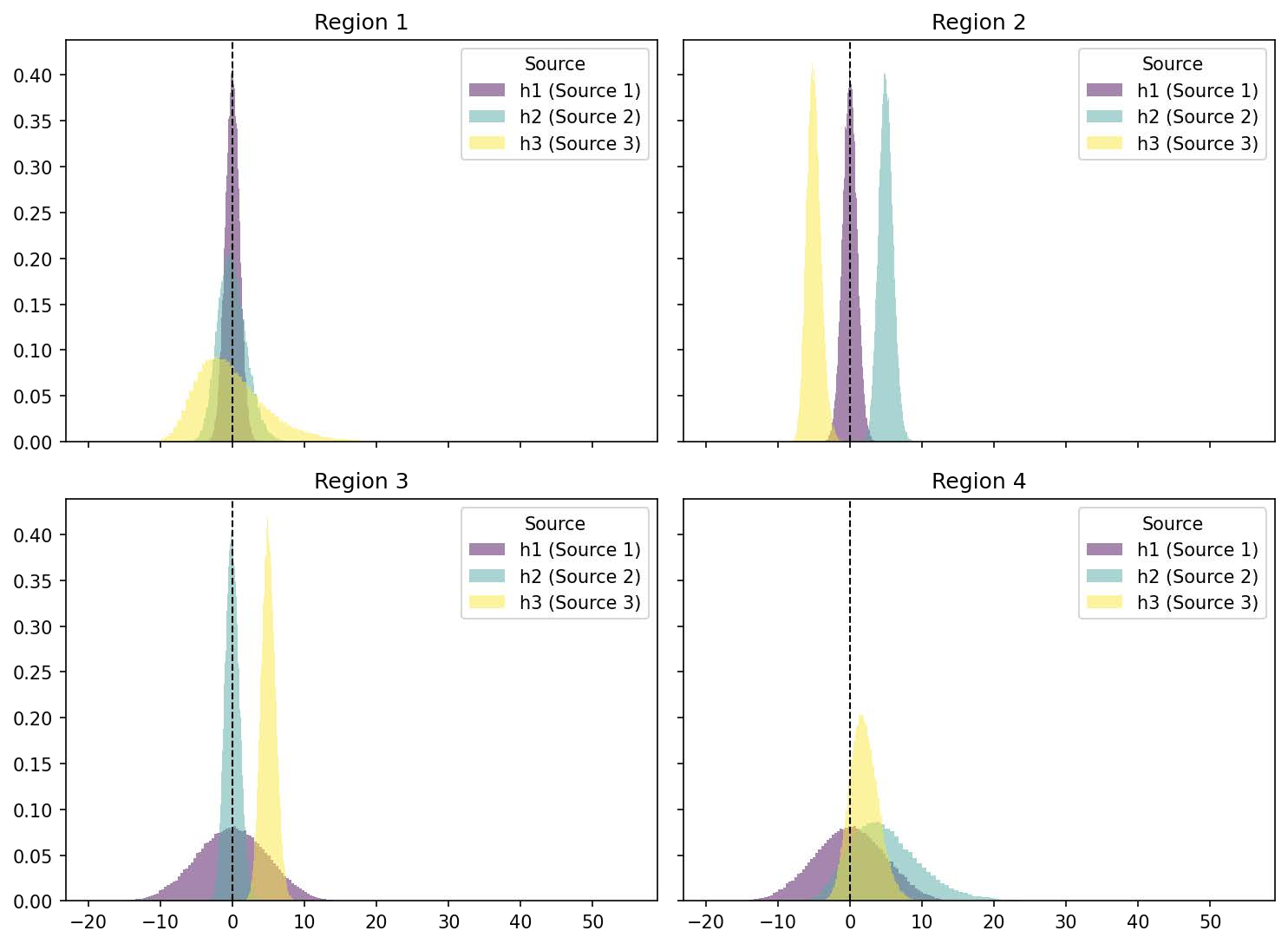}
    \caption{Prediction error distributions for the resource allocation problem (varying distribution types)}
    \label{fig:error-distribution-type}
\end{figure}

    One of the advantages of our approach that constructs the ambiguity set in this study is that it does not require that the prediction errors from different information sources follow the same type of probability distribution. To test this, we modify the baseline setting such that the prediction values provided by Source $2$ and Source $3$ are sampled from Lognormal distributions instead of truncated Normal distributions (see Figure~\ref{fig:error-distribution-type}). The results, reported in Table~\ref{tab:performance-resource-distribution}, show that varying the distribution types does not influence the performance of MR-DRO, and it still outperforms other benchmark models that rely on a single information source. 

\begin{table}[htbp!]
\centering
\caption{Performances of different methods for the resource allocation problem with varying distribution types}
\begin{tabular}{c c c c c} 
\toprule
Model & Obj. (*\$1000) &  Loss (*\$1000) & Time (sec.)\\
\midrule
MR-DRO (Min-max) & $7.725 \pm 0.453$ & $8.259 \pm 0.106$ & $129.874 \pm 34.725$\\
MR-DRO (Exponential) & $6.942 \pm 0.148$ & $7.099 \pm 0.101$ & $116.901 \pm 35.571$\\
MR-DRO (Variable-share) & $8.732 \pm 1.360$ & $7.207 \pm 0.116$ & $135.199 \pm 35.179$\\
DRO (h1) & $13.801 \pm 2.058$ & $15.653 \pm 0.216$ & $111.402 \pm 35.087$\\
DRO (h2) & $12.287 \pm 0.553$ & $12.762 \pm 0.093$ & $110.515 \pm 35.257$\\
DRO (h3) & $10.496 \pm 0.859$ & $11.260 \pm 0.080$ & $111.331 \pm 35.229$\\
\bottomrule
\end{tabular}
\label{tab:performance-resource-distribution}
\end{table}

\paragraph{Nonstationary error distributions.}

    Another advantage of trust dynamics is its ability to adapt when the error distribution of each source is nonstationary, changing throughout the trust update process. In this experiment, each source's error distribution shifts at predefined time points, reflecting real-world scenarios where the accuracy of a reference information source can fluctuate over time. We set the original mean and standard deviation of the error distribution of each source as $\bm{\mu} = (\bm{\mu}_{1},\bm{\mu}_{2},\bm{\mu}_{3})^\mathsf{T} = ((0,0,0,0)^\mathsf{T},(0,5,0,5)^\mathsf{T},(0,-5,5,2)^\mathsf{T})^\mathsf{T}$ and $\bm{\sigma} = (\bm{\sigma}_{1},\bm{\sigma}_{2},\bm{\sigma}_{3})^\mathsf{T} = ((1,1, 5,5)^\mathsf{T},(2,1,1,5)^\mathsf{T},(5,1,1,2)^\mathsf{T})^\mathsf{T}$. The predefined time points are $(100, 100, 100, 50)$ in each region. After the corresponding time points, $\bm{\mu}$ and $\bm{\sigma}$ change to $\bm{\mu}^{'} = (\bm{\mu}_{1}^{'},\bm{\mu}_{2}^{'},\bm{\mu}_{3}^{'})^\mathsf{T} = ((5,0,0,5)^\mathsf{T},(0,0,0,0)^\mathsf{T},(0,-5,0,0)^\mathsf{T})^\mathsf{T}$ and $\bm{\sigma}^{'} = (\bm{\sigma}_{1}^{'},\bm{\sigma}_{2}^{'},\bm{\sigma}_{3}^{'})^\mathsf{T} = ((1,1,5,1)^\mathsf{T},(2,1,2,2)^\mathsf{T},(5,1,1,5)^\mathsf{T})^\mathsf{T}$. For example, in Region $1$, the error distribution of source $1$ shifts from a Normal distribution $\mathcal{N}(0,1)$ to $\mathcal{N}(5,1)$ after $100$ events have occurred.

\begin{figure}[ht]
    \centering
    \begin{subfigure}{\textwidth}
        \centering
        \includegraphics[width=\textwidth]{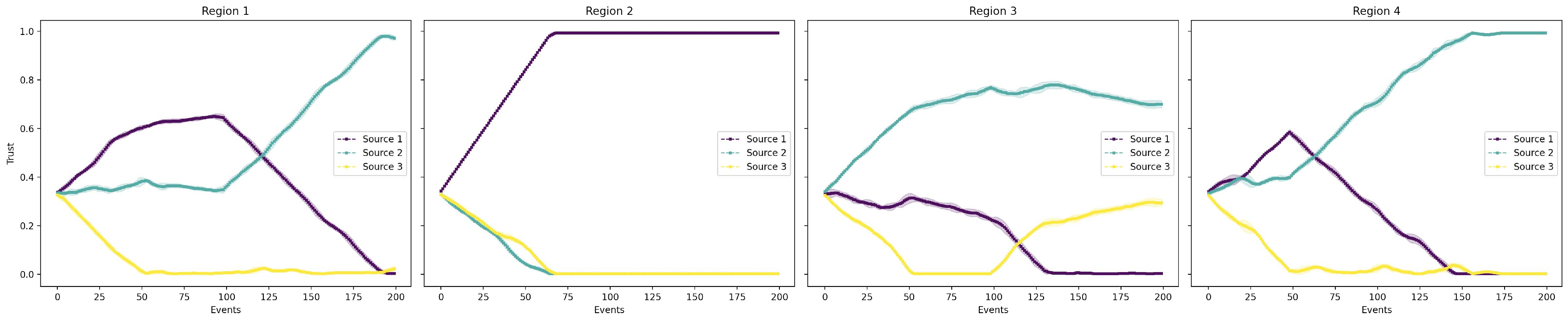}
        \caption{Min-max error trust update}
        \label{fig:resource_trust_update_min_max_nonstationary}
    \end{subfigure}   
    \vspace{0.5em} 
    \begin{subfigure}{\textwidth}
        \centering
        \includegraphics[width=\textwidth]{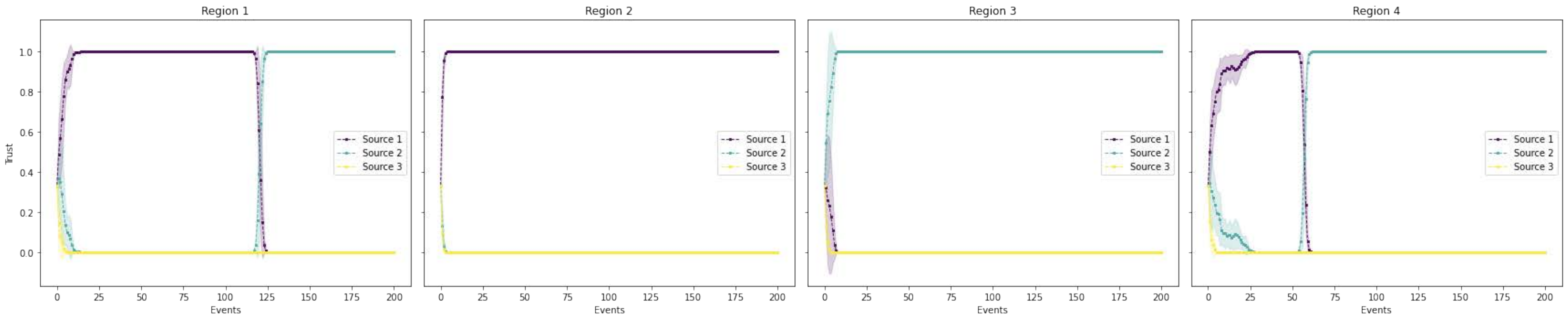}
        \caption{Exponential error trust update}
        \label{fig:resource_trust_update_exponential_nonstationary}
    \end{subfigure}   
    \vspace{0.5em} 
    \begin{subfigure}{\textwidth}
        \centering
        \includegraphics[width=\textwidth]{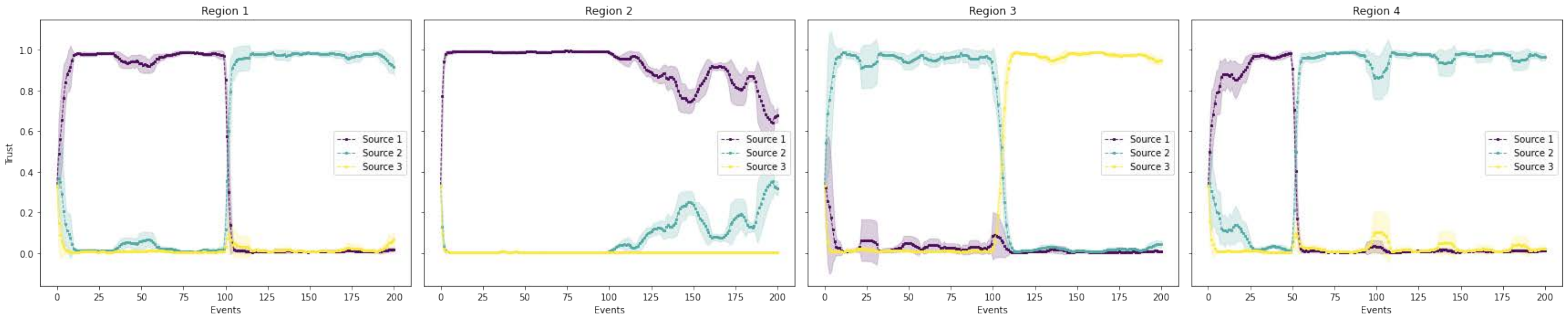}
        \caption{Variable-share error trust update}
        \label{fig:resource_trust_update_variable_share_nonstationary}
    \end{subfigure}
    \caption{Trust update process in the nonstationary error distribution setting for the resource allocation problem ($K_{\text{r}} = 4$, $H = 3$, $I = 200$).}
    \label{fig:resource_trust_update__nonstationary}
\end{figure}

    The following results of the experiment demonstrate how the trust dynamics responds to these changes, adjusting trust levels based on the changing performance of each source, as illustrated in Figure~\ref{fig:resource_trust_update__nonstationary}. The figure highlights how trust dynamics effectively down-weight or up-weight sources as their error profiles change, thus maintaining robust decision-making even when source reliability is not consistent over time. For example, before the predefined time point, source $1$ is the source with the highest accuracy rate among all sources in region $1$. However, after the time point, the prediction error mean of Source $1$ increases, making it an unreliable reference source compared to Source $2$. Thus, the trust in Source $1$ decreases after its accuracy rate decreases. Table~\ref{tab:performance-resource-nonstationary-distribution} summarizes the performance of different methods. In particular, in Region $2$ and Region $3$, the variable-share error update algorithm manages to detect the prediction accuracy shift in sources better than the min-max and exponential error update algorithm, which aligns with the design purpose of this algorithm.

\begin{table}[htbp!]
\centering
\caption{Performances of different methods for the resource allocation problem with non-stationary error distributions}
\begin{tabular}{c c c c c} 
\toprule
Model & Obj. (*\$1000) &  Loss (*\$1000) & Time (sec.)\\
\midrule
MR-DRO (Min-max) & $10.918 \pm 0.384$ & $10.820 \pm 0.069$ & $128.459 \pm 31.378$\\
MR-DRO (Exponential) & $6.942 \pm 0.148$ & $16.122 \pm 0.101$ & $116.901 \pm 35.571$\\
MR-DRO (Variable-share) & $8.732 \pm 1.360$ & $7.207 \pm 0.116$ & $135.199 \pm 35.179$\\
DRO (h1) & $15.602 \pm 0.460$ & $16.122 \pm 0.126$ & $103.481 \pm 31.254$\\
DRO (h2) & $11.012 \pm 1.084$ & $11.331 \pm 0.069$ & $103.446 \pm 31.388$\\
DRO (h3) & $16.416 \pm 1.940$ & $19.129 \pm 0.159$ & $102.966 \pm 31.380$\\
\bottomrule
\end{tabular}
\label{tab:performance-resource-nonstationary-distribution}
\end{table}

\subsection{Probability dominance and dominant trust}\label{sec:5.2}
In this section, we numerically validate the result of probability dominance in Theorem~\ref{thm:theorem4} and illustrate the related trust update patterns, by focusing on special cases with $H=2$ information sources. 
For $K_r =4$ regions, we set $\bm{\mu} = (\bm{\mu}_{1},\bm{\mu}_{2})^\mathsf{T} = ((0,0,0,2)^\mathsf{T},(0,5,2,-2)^\mathsf{T})^\mathsf{T}$ and $\bm{\sigma} = (\bm{\sigma}_{1},\bm{\sigma}_{2})^\mathsf{T} = ((1,5,5,2)^\mathsf{T},(5,5,2,2)^\mathsf{T})^\mathsf{T}$. We observe the trust update process for $I=300$ events. Judged from the empirical cumulative distribution function and based on Theorem \ref{thm:theorem4}, the deviation from source $1$ dominates the deviation from source $2$ with probability $\beta=0.5$; while in region $3$, the deviation from source $2$ has probability dominance with $\beta=0.5$ over the deviation from source $1$. In region $4$, there is no clear probability dominance between these two sources.

\begin{figure}[htbp!]
    \centering
    \includegraphics[width=\textwidth]{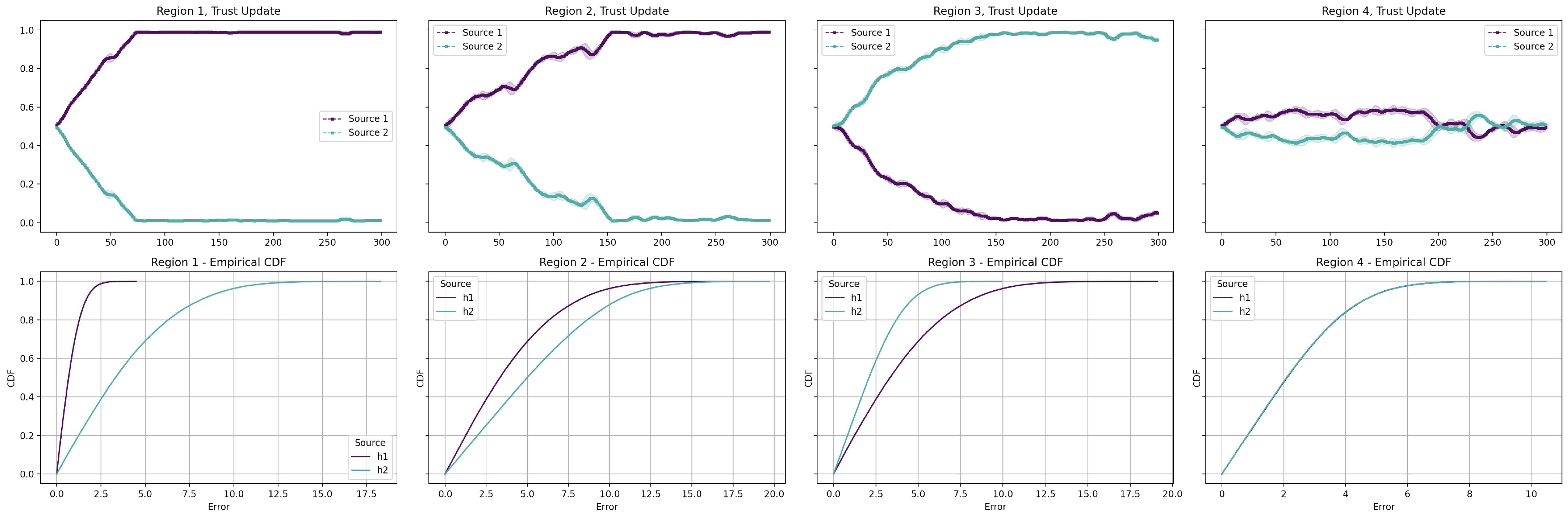}
    \caption{Trust update process with/without probability dominance relationship between two sources}
    \label{fig:fsd}
\end{figure}

Figure~\ref{fig:fsd} depicts the trust update results for the $30$ trials. In regions where the trust for one source dominates the other, there exists a corresponding probability dominance relationship between the sources. In contrast, in regions without a clear probability dominance relationship, we cannot observe a dominant trust pattern along the trust update process. These results are aligned with Theorem~\ref{thm:theorem4}. 

\subsection{Portfolio optimization with uncertain return}\label{sec:5.4}
We provide another set of numerical studies based on portfolio optimization instances, where the uncertainty appears in the objective function. 

\subsubsection{MR-DRO for a portfolio optimization problem}\label{sec:5.4.1}
Consider a capital market with $K_{\text{a}}$ assets, where the returns of these assets are characterized by the random vector $\bm{r} = [r_{1},\ldots,r_{K_{\text{a}}}]^\mathsf{T}$, $\bm{r} \in \mathbb{R}^{K_{\text{a}}}$. Short-selling forbidden, a portfolio can be encoded by a vector $\bm{x}=[x_{1},\ldots, x_{K_{\text{a}}}]^\mathsf{T}$ of percentage weights of each asset. For each asset $k \in [K_{\text{a}}]$, a proportion $x_{k}$ of the total capital is invested, leading to a portfolio return of $\langle\bm{x},\bm{r}\rangle$. Accordingly, we solve the  the following MR-DRO model:
\begin{equation}   
    J_{\textrm{portfolio}}^{*} = \inf_{\bm{x} \in \mathbb{X}}\left\{\sup_{\mathbb{Q} \in \mathbb{B}_{\epsilon}(\hat{\mathbb{P}}_{HI})}\left\{\mathbb{E}^{\mathbb{Q}}\left[-\langle\bm{x},\bm{r}\rangle\right] + \rho \mathbb{Q}-\textrm{CVaR}_{\alpha}(-\langle\bm{x},\bm{r}\rangle)\right\}\right\},
\end{equation}
which minimizes a weighted sum of the expectation and the conditional value at-risk (CVaR) of the portfolio loss $-\langle\bm{x},\bm{r}\rangle$, where $\alpha \in (0,1]$ is referred to as the confidence level of the CVaR, and $\rho \in \mathbb{R}_{+}$ quantifies the investor's risk-aversion. The CVaR at the $\alpha$ level essentially captures the average of the worst $\alpha \times 100\%$ portfolio losses under distribution $\mathbb{Q}$. By expressing the CVaR with its formal definition, we arrive at
\begin{align}
    J_{\textrm{portfolio}}^{*} &= \inf_{\bm{x} \in \mathbb{X}}\left\{\sup_{\mathbb{Q} \in \mathbb{B}_{\epsilon}(\hat{\mathbb{P}}_{HI})}\left\{\mathbb{E}^{\mathbb{Q}}\left[-\langle\bm{x},\bm{r}\rangle\right] + \rho \inf_{\tau \in \mathbb{R}} \left\{\mathbb{E}^{\mathbb{Q}}\left[\tau + \frac{1}{\alpha}\max\{-\langle\bm{x},\bm{r}\rangle-\tau,0\}\right]\right\}\right\}\right\}\nonumber\\
    &= \inf_{\bm{x} \in \mathbb{X}, \tau \in \mathbb{R}}\left\{\sup_{\mathbb{Q} \in \mathbb{B}_{\epsilon}(\hat{\mathbb{P}}_{HI})}\left\{\mathbb{E}^{\mathbb{Q}}\left[\max_{j \in [J]}(\langle\bm{a}_{j},\bm{r}\rangle + b_{j})\right]\right\}\right\},
\end{align}
where $J=2$, $\bm{a}_{1} = -\bm{x}$, $a_{2} = (-1-\rho/\alpha)\bm{x}$, $b_{1} = \rho\tau$, $b_{2}=\rho\tau(1-1/\alpha)$. In a scenario where an investor does not have access to the true distribution $\mathbb{P}$ but possesses the empirical distribution $\hat{\mathbb{P}}_{HI}$ constructed by nonparametric data fusion based on trust, we can apply Theorem~\ref{thm:piecewise-affine} to determine an optimal portfolio $\bm{x}$ and formulate the problem as \eqref{eq:DRO-portfolio}:
\begin{subequations}\label{eq:DRO-portfolio}
\begin{align}
    \inf_{\bm{x}, \lambda, s_{hi}, \bm{\gamma}_{hij}}\quad & \lambda\epsilon + \sum_{h=1}^{H}\sum_{i=1}^{I}\frac{t^{(I)}_{h}}{I}s_{hi}\label{eq:DRO-portfolio-obj}&\\
    \textrm{s.t.}\quad\quad\quad & \sum_{k=1}^{K}x_{k}=1,\label{eq:DRO-portfolio-x}&\\
    & b_{j} + a_{j}\langle\bm{x},\hat{\bm{r}}^{(i)}_{h}\rangle, + \langle \gamma_{hij}, d - C\hat{\bm{r}}^{(i)}_{h} \rangle \leq s_{hi}, &\forall i \in [I], j \in [J], h \in [H],\label{eq:DRO-portfolio-constr1}\\
    & \|C^\mathsf{T}\bm{\gamma}_{hij} - a_{j}\bm{x}\|_{*} \leq \lambda, &\forall i \in [I], j \in [J], h \in [H],\label{eq:DRO-portfolio-constr2}\\
    & x_{k} \geq0, \bm{\gamma}_{hij} \geq 0, &\forall i \in [I], j \in [J], k \in [K_{\text{a}}], h \in [H].\label{eq:DRO-portfolio-constr3}
\end{align}
\end{subequations}
\subsubsection{Baseline settings for the portfolio optimization problem}\label{sec:5.4.2}
We use NASDAQ100 dataset in \citet{bruni2016real} with $K_{\text{a}}=82$. The true return $\bm{r}^{(i)}_{\textrm{true}}$ at time period $i$ is a vector with each element equals the true return value corresponding to the dataset. The predicted return of asset $k$ from source $h$ at time period $i$ is sampled from a truncated normal distribution with mean $(r^{(i)}_{\textrm{true},k}+\mu_{h,k})$ and variance $\sigma^{2}_{h,k}$ within interval $(-1,1)$, which means no asset can lose more than $100\%$ of its value. We set $\alpha = 20\%$ and $\rho=10$ in all numerical experiments. The step size for min-max error trust update is $\Delta t=0.01$. The update rate is $\eta=10^{2}$ for exponential error trust update as well as variable-share error trust update, and the share portion is $\beta = 0.5$ in variable-share trust update. We conduct $30$ trials with distinct random seeds for each trust update algorithm to observe common patterns during the trust update process.

\subsubsection{Baseline results for the portfolio optimization problem}\label{sec:5.4.3}

\begin{figure}[htbp!]
    \centering
    \includegraphics[width=0.85\textwidth]{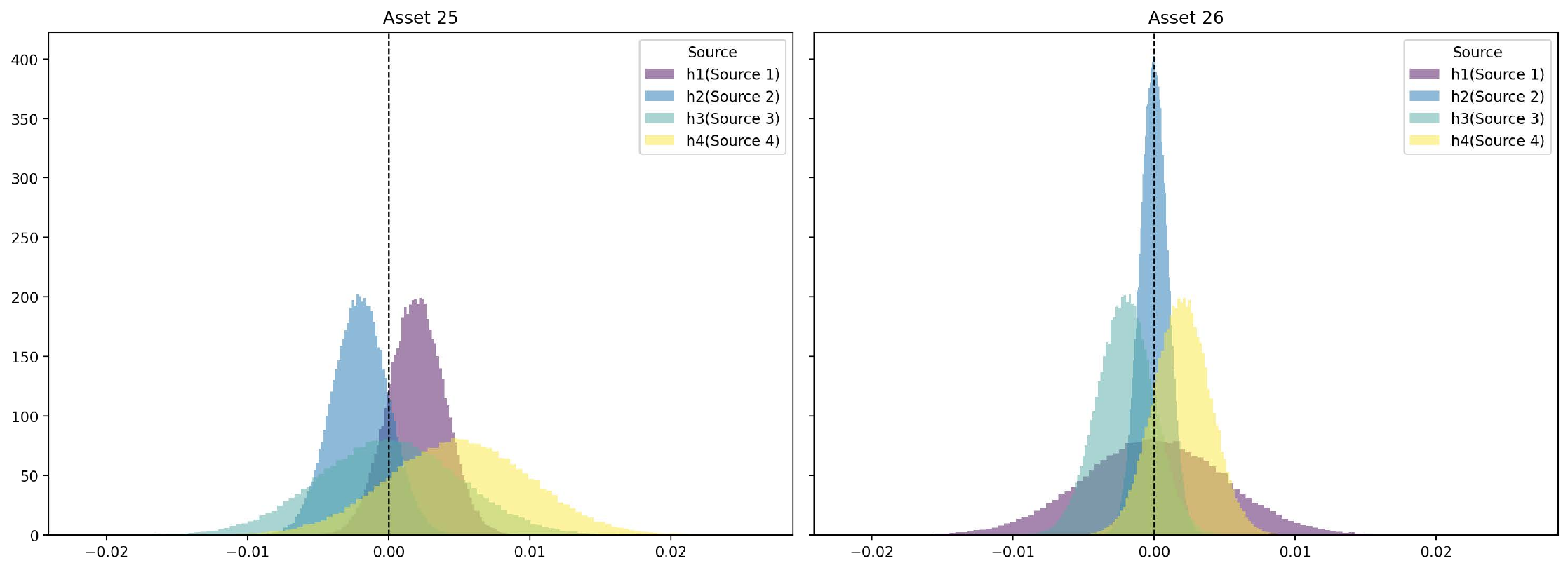}
    \caption{Prediction error distribution with the baseline setting for the portfolio optimization problem}
    \label{fig:error-distribution-portfolio}
\end{figure}

\begin{figure}[htbp!]
    \centering
    \begin{subfigure}{\textwidth}
        \centering
        \includegraphics[width=\textwidth]{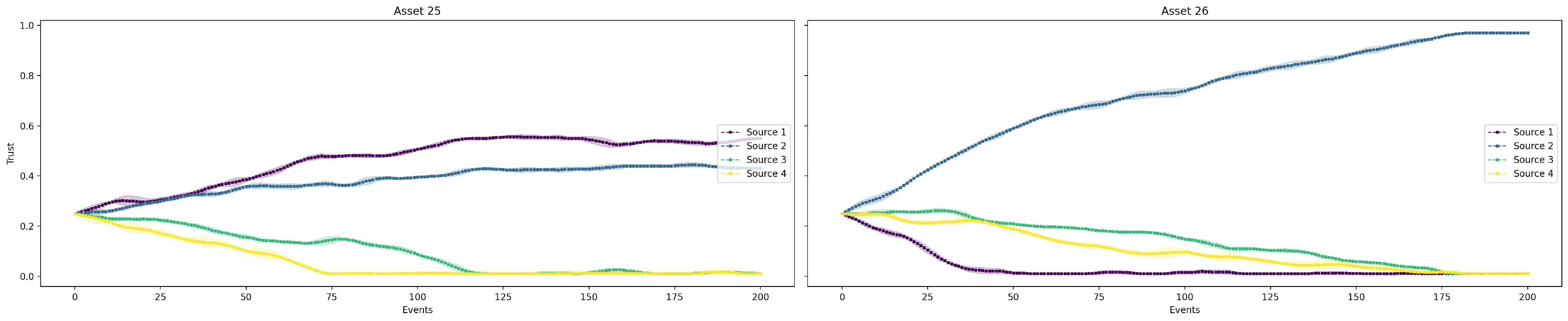}
        \caption{Min-max error trust update}
        \label{fig:portfolio_trust_update_min_max}
    \end{subfigure} 
    \vspace{0.5em} 
    \begin{subfigure}{\textwidth}
        \centering
        \includegraphics[width=\textwidth]{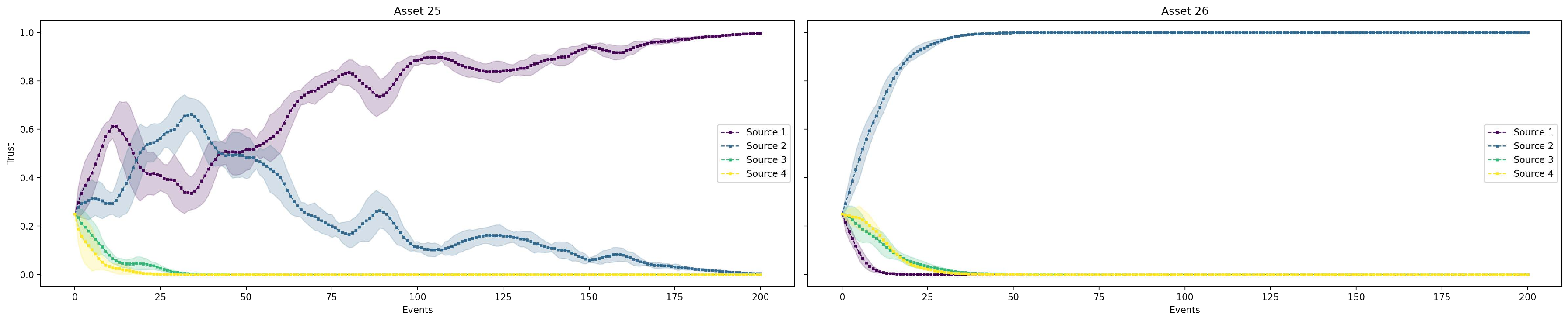}
        \caption{Exponential error trust update}
        \label{fig:portfolio_trust_update_exponential}
    \end{subfigure}
    \vspace{0.5em} 
    \begin{subfigure}{\textwidth}
        \centering
        \includegraphics[width=\textwidth]{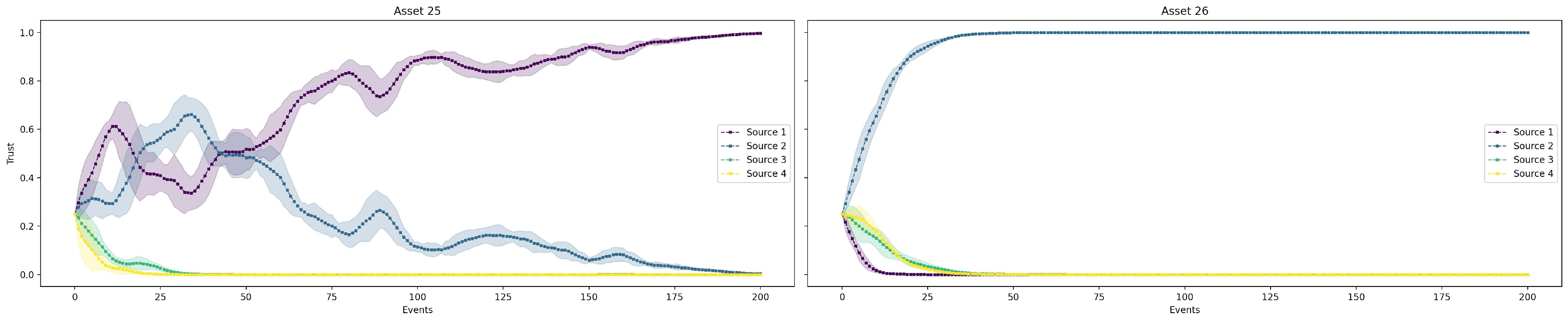}
        \caption{Variable-share error trust update}
        \label{fig:portfolio_trust_update_variable_share}
    \end{subfigure}  
    \caption{Trust Update Process in the baseline setting for the portfolio optimization problem ($K_{\text{a}} = 28$, $H = 4$, $I = 200$).}
    \label{fig:DowJones}
\end{figure}

Using the baseline setting and Algorithms~\ref{alg:min-max-error-trust-update}--\ref{alg:variable-share-error-trust-update}, we obtain the result of the trust update process, reported in Table~\ref{tab:performance-portfolio}. Due to the problem scale, we only include the error distribution demonstration and trust update results for Asset $25$ and $26$ in Figure~\ref{fig:error-distribution-portfolio} and Figure~\ref{fig:DowJones}, as they are typical examples of the absence or presence of probabilistic dominance relationships between information sources, which will influence whether or not the dominant trust pattern will show when using the min-max error trust update algorithm.

\begin{table}[htbp!]
\centering
\caption{Performances of different methods with the baseline setting for the portfolio optimization problem ($K_{\text{a}} = 28$, $H=4$ and $I=200$)}
\begin{tabular}{c c c c} 
\toprule
 Model & Obj. & Loss & Time (sec.)\\
\midrule
 MR-DRO(Min-max) & $0.018\pm8.223 \times 10^{-6}$ & $0.002\pm0.233 \times 10^{-3}$ & $472.177 \pm 5.608$\\
 MR-DRO(Exponential) & $0.018\pm8.225 \times 10^{-6}$ & $0.002\pm0.218\times 10^{-3}$ & $444.359\pm7.746$\\
 MR-DRO(Variable-share) & $0.018\pm8.224 \times 10^{-6}$ & $0.002\pm0.219 \times 10^{-3}$ & $489.041\pm6.029$\\
 DRO (h1) & $0.018\pm1.324 \times 10^{-5}$ & $0.004\pm0.703 \times 10^{-3}$ & $247.786 \pm 3.377$\\
 DRO (h2) & $0.018\pm8.312 \times 10^{-6}$ & $0.004\pm0.799 \times 10^{-3}$ & $251.808 \pm 3.679$\\
 DRO (h3) & $0.018\pm1.001 \times 10^{-5}$ & $0.005\pm0.673 \times 10^{-3}$ & $248.473 \pm 3.262$\\
 DRO (h4) & $0.018\pm1.071 \times 10^{-5}$ & $0.004\pm0.635 \times 10^{-3}$ & $248.997 \pm 3.396$\\
\bottomrule
\end{tabular}
\label{tab:performance-portfolio}
\end{table}

In this setting, the MR-DRO models still perform better than using DRO with single-source data in both the objective value and the average loss, with a slightly longer computation time.

We also performed out-of-sample performance analysis, where for given solution $\bm{x}^{*}$, the performance metric is given by: 
\begin{equation}
    J_{\textrm{portfolio}}(\bm{x}^{*}) :=\frac{1}{|\Psi|}\sum_{\psi \in \Psi}\left[-\langle\bm{x}^{*},\bm{r}^{\psi}\rangle\right] + \rho \times \textrm{CVaR}_{\alpha}\left(\left\{-\langle\bm{x}^{*},\bm{r}^{\psi}\rangle\right\}_{\psi \in \Psi}\right).\nonumber
\end{equation}
MR-DRO models consistently have the best performance among all methods, reported in Table~\ref{tab:out-of-sample-performance-portfolio}.

\begin{table}[htbp!]
\centering
\caption{Out-of-sample performances of different methods for the portfolio optimization problem ($K_{\text{a}}=28$ and $H=4$)}
\begin{tabular}{c c c} 
\toprule
 $|\Psi|$ & Model & Loss\\
\midrule
 \multirow{ 7}{*}{$40$} & MR-DRO (Min-max) & $-0.031\pm0.491 \times 10^{-3}$\\
  & MR-DRO (Exponential) & $-0.034\pm0.545 \times 10^{-3}$\\
  & MR-DRO (Variable-share) & $-0.034\pm0.519 \times 10^{-3}$\\
  & DRO (h1) & $-0.012\pm1.052 \times 10^{-3}$\\
  & DRO (h2) & $-0.012\pm1.517 \times 10^{-3}$\\
  & DRO (h3) & $-0.008\pm1.627 \times 10^{-3}$\\
  & DRO (h4) & $-0.025\pm1.070 \times 10^{-3}$\\
\bottomrule
\end{tabular}
\label{tab:out-of-sample-performance-portfolio}
\end{table}

\subsection{Comparative insights from numerical studies}\label{sec:5.5}
We compare the results of experiments carried out in two distinct contexts: resource allocation and portfolio optimization. Both sets of experiments demonstrate the effectiveness of MR-DRO models, irrespective of the placement of uncertainty within the mathematical model. In the resource allocation problem, uncertainty is embedded originally within the constraints, specifically in the demand parameters, while in the portfolio optimization problem, it resides in the objective function, influencing the returns of assets. Despite these structural differences, the MR-DRO models consistently outperform the DRO models with reference information from a single source. This indicates a high degree of flexibility and adaptability in different settings, as MR-DRO models adapt to the specific characteristics of each problem. Furthermore, when the pairwise prediction error distribution relationship satisfies Theorem~\ref{thm:theorem5} in one region or asset of the problem, we observe the emergence of dominant trust in both contexts when using the min-max error trust update algorithm.

\section{Conclusion}\label{sec:conclusion}
In this paper, we developed and promoted the use of the MR-DRO approach to tackle optimization under uncertainty having parameter information from multiple heterogeneous data sources. Our model features a nonparametric data fusion Wasserstein ambiguity set, which improves decision accuracy by reducing the impact of errors from different information sources. We proposed a trust update mechanism with three different trust update algorithms based on historical errors and introduced the concept of probability dominance to explain the decision maker's preference for one source over another. We demonstrated our model's effectiveness through its application in a resource allocation and a portfolio optimization problem, where it consistently outperformed DRO models that rely on a single information source, particularly in minimizing worst-case losses. To the best of our knowledge, this is one of the very first endeavors to bring the notion of trust in human-computer interaction to decision models and frameworks for optimization under uncertainty.


\ACKNOWLEDGMENT{This work is supported in part by the United States (US) Army Research Laboratory under contract no. W911NF2320215 and the US National Science Foundation (NSF) grant no. ECCS-2533775.}

\bibliographystyle{apacite}
\bibliography{MR-DRO.bib} 

\newpage
\setcounter{page}{1}
\begin{appendices}

\section{Supplement Materials}\label{sec:proof}
\subsection{Proof of Theorem~\ref{thm:convex-reduction}}\label{sec:proof-1}
\proof{Proof:}
Following~\eqref{eq:Wasserstein-metric} we rewrite the worst case expectation \eqref{eq:worst-case-expectation} as 
\begin{align}
    \sup_{\mathbb{Q} \in \mathbb{B}_{\epsilon}(\hat{\mathbb{P}}_{HI})}\mathbb{E}^{\mathbb{Q}}\left[\ell(\bm{x},\bm{\xi})\right] &=
    \begin{cases}
    \sup_{\Pi, \mathbb{Q}} \int_{\Xi} \ell(\bm{x},\bm{\xi}) \mathbb{Q}(d\bm{\xi})\\
    \textrm{s.t.} \quad \int_{\Xi^{2}} \|\bm{\xi} - \bm{\xi}^{'}\| \Pi(d\bm{\xi},d\bm{\xi}^{'}) \leq \epsilon
    \end{cases}\label{eq:first-equality}
\\&=
    \begin{cases}
    \sup_{\mathbb{Q}^{(i)}_{h} \in \mathcal{M}(\Xi)} \sum_{h=1}^{H}\sum_{i=1}^{I}\frac{t^{(I)}_{h}}{I}\int_{\Xi} \ell(\bm{x},\bm{\xi}) \mathbb{Q}^{(i)}_{h}(d\bm{\xi}) \\
    \textrm{s.t.} \quad \sum_{h=1}^{H}\sum_{i=1}^{I}\frac{t^{(I)}_{h}}{I}\int_{\Xi} \|\bm{\xi}-\hat{\bm{\xi}}^{(i)}_{h}\| \mathbb{Q}^{(i)}_{h}(d\bm{\xi}) \leq \epsilon,           
    \end{cases}\label{eq:second-equality}
\end{align}
where $\Pi$ is a joint distribution of $\bm{\xi}$ and $\bm{\xi}^{'}$ with marginals $\mathbb{Q}$ and $\hat{\mathbb{P}}_{HI}$, respectively. We drop the minimization problem in \eqref{eq:Wasserstein-metric} and obtain \eqref{eq:first-equality} since the minimization of $d_{W}(\mathbb{Q},\hat{\mathbb{P}}_{HI})$ is at most radius $\epsilon$ equivalent as \eqref{eq:first-equality} has a feasible solution. The second equality \eqref{eq:second-equality} follows from the law of total probability, which means that any joint probability distribution $\Pi$ of $\bm{\xi}$ and $\bm{\xi}^{'}$ can be constructed from the marginal distribution $\hat{\mathbb{P}}_{HI}$ of $\bm{\xi}^{'}$ and the conditional distributions $\mathbb{Q}^{(i)}_{h}$ of $\bm{\xi}$ given $\bm{\xi}^{'} = \hat{\bm{\xi}}^{(i)}_{h}$, for all $i \in [I]$ and $h \in [H]$. In other words, we have $\Pi = \sum_{h=1}^{H}\sum_{i=1}^{I} \frac{t^{(I)}_{h}}{I}\cdot \delta(\bm{\xi} - \hat{\bm{\xi}}^{(i)}_{h}) \bigotimes \mathbb{Q}^{(i)}_{h}$, where $\bigotimes$ denotes the Kronecker product.

Applying the standard duality argument \citep{bertsimas1997introduction}, we obtain
\begin{subequations}
\begin{align}
\sup_{\mathbb{Q} \in \mathbb{B}_{\epsilon}(\hat{\mathbb{P}}_{HI})}\mathbb{E}^{\mathbb{Q}}\left[\ell(\bm{x},\bm{\xi})\right] &=
\sup_{\mathbb{Q}^{(i)}_{h} \in \mathcal{M}(\Xi)} \inf_{\lambda \geq 0} \Bigg\{\sum_{h=1}^{H}\sum_{i=1}^{I}\frac{t^{(I)}_{h}}{I} \int_{\Xi} \ell(\bm{x},\bm{\xi}) \mathbb{Q}^{(i)}_{h}(d\bm{\xi})\nonumber \\
& \quad \quad \quad \quad \quad \quad + \lambda \left(\epsilon - \sum_{h=1}^{H}\sum_{i=1}^{I} \frac{t^{(I)}_{h}}{I} \int_{\Xi} \|\bm{\xi}-\hat{\bm{\xi}}^{(i)}_{h}\| \mathbb{Q}^{(i)}_{h}(d\bm{\xi})\right) \Bigg\} \nonumber \\
& \leq
\inf_{\lambda \geq 0} \sup_{\mathbb{Q}^{(i)}_{h} \in \mathcal{M}(\Xi)} \left\{ \lambda\epsilon + \sum_{h=1}^{H}\sum_{i=1}^{I}\frac{t^{(I)}_{h}}{I}\int_{\Xi}\left(\ell(\bm{x},\bm{\xi}) - \lambda\|\bm{\xi}-\hat{\bm{\xi}}^{(i)}_{h}\|\right) \mathbb{Q}^{(i)}_{h}(d\bm{\xi})\right\} \label{eq:max-min} \\
& =
\inf_{\lambda \geq 0} \quad \left\{\lambda\epsilon + \sum_{h=1}^{H}\sum_{i=1}^{I} \frac{t^{(I)}_{h}}{I} \sup_{\bm{\xi} \in \Xi} \left\{\ell(\bm{x},\bm{\xi}) - \lambda\|\bm{\xi}-\hat{\bm{\xi}}^{(i)}_{h}\|\right\}\right\},\label{eq:Dirac}
\end{align}
\end{subequations}
where \eqref{eq:max-min} holds because of the max-min inequality, and \eqref{eq:Dirac} holds as $\mathcal{M}(\Xi)$ contains all the Dirac distributions supported on $\Xi$. After introducing epigraphical auxiliary variables $s_{hi}$, $h \in [H]$, $i \in [I]$, we reformulate \eqref{eq:Dirac} as:
\begin{subequations}
\begin{align}
    & \quad
    \begin{cases}
        \inf_{\lambda,s_{hi}} \quad \lambda\epsilon + \sum_{h=1}^{H}\sum_{i=1}^{I}\frac{t^{(I)}_{h}}{I}s_{hi}\\
        \textrm{s.t.} \quad \sup_{\bm{\xi} \in \Xi}\left\{\ell_{j}(\bm{x},\bm{\xi})-\lambda\|\bm{\xi} - \hat{\bm{\xi}}^{(i)}_{h}\|\right\} \leq s_{hi}, \quad \forall h \in [H], i \in [I], j \in [J],\\
        \quad \quad \lambda \geq 0
    \end{cases}\label{eq:introduce-s}\\
    & =
    \begin{cases}
        \inf_{\lambda,s_{hi},\bm{z}_{hij}} \quad \lambda\epsilon + \sum_{h=1}^{H}\sum_{i=1}^{I}\frac{t^{(I)}_{h}}{I}s_{hi}\\
        \textrm{s.t.} \quad \sup_{\bm{\xi} \in \Xi}\left\{\ell_{j}(\bm{x},\bm{\xi}) - \max_{\|\bm{z}_{hij}\|_{*} \leq \lambda}\langle \bm{z}_{hij},\bm{\xi} - \hat{\bm{\xi}}^{(i)}_{h}\rangle\right\} \leq s_{hi},\quad \forall h \in [H], i \in [I], j \in [J],\\
        \quad \quad \lambda \geq 0
    \end{cases}\label{eq:dual-norm-definition}\\
    & \leq
    \begin{cases}
        \inf_{\lambda,s_{hi},\bm{z}_{hij}} \quad \lambda\epsilon + \sum_{h=1}^{H}\sum_{i=1}^{I}\frac{t^{(I)}_{h}}{I}s_{hi}\\
        \textrm{s.t.} \quad \sup_{\bm{\xi} \in \Xi}\left\{\ell_{j}(\bm{x},\bm{\xi})-\langle \bm{z}_{hij},\bm{\xi}\rangle\right\}+\langle \bm{z}_{hij},\hat{\bm{\xi}}^{(i)}_{h}\rangle \leq s_{hi}, \quad \forall h \in [H], i \in [I], j \in [J],\\
        \quad \quad \|\bm{z}_{hij}\|_{*} \leq \lambda, \quad \forall h \in [H], i \in [I], j \in [J],
    \end{cases}\label{eq:rewrite-interchange}\\
    & =
    \begin{cases}
        \inf_{\lambda,s_{hi},\bm{z}_{hij}} \quad \lambda\epsilon + \sum_{h=1}^{H}\sum_{i=1}^{I}\frac{t^{(I)}_{h}}{I}s_{hi}\\
        \textrm{s.t.} \quad [-\ell_{j} + \chi_{\Xi}]^{*}(-\bm{z}_{hij})+\langle \bm{z}_{hij},\hat{\bm{\xi}}^{(i)}_{h}\rangle \leq s_{hi}, \quad \forall h \in [H], i \in [I], j \in [J],\\
        \quad \quad \|\bm{z}_{hij}\|_{*} \leq \lambda, \quad \forall h \in [H], i \in [I], j \in [J],   
    \end{cases}\label{eq:conjugacy}\\
    & = 
    \begin{cases}
        \inf_{\lambda,s_{hi},\bm{z}_{hij}} \quad \lambda\epsilon + \sum_{h=1}^{H}\sum_{i=1}^{I}\frac{t^{(I)}_{h}}{I}s_{hi}\\
        \textrm{s.t.} \quad [-\ell_{j} + \chi_{\Xi}]^{*}(\bm{z}_{hij})-\langle \bm{z}_{hij},\hat{\bm{\xi}}^{(i)}_{h}\rangle \leq s_{hi}, \quad \forall h \in [H], i \in [I], j \in [J],\\
        \quad \quad \|\bm{z}_{hij}\|_{*} \leq \lambda, \quad \forall h \in [H], i \in [I], j \in [J],\label{eq:substitue-minus-z}\\    
    \end{cases}
\end{align}
\end{subequations}
We obtain \eqref{eq:dual-norm-definition} from \eqref{eq:introduce-s} by the definition of the dual norm $\|\bm{z}_{hij}\|_{*}:=\sup_{\|\bm{u}\| \leq 1}\langle \bm{z}_{hij},\bm{u}\rangle$. This measures the maximum projection of 
$\bm{z}_{hij} \in \mathbb{R}^{M}$ onto any unit vector $\bm{u} \in \mathbb{R}^{M}$. By introducing the constraint $\|\bm{z}_{hij}\|_{*} \leq \lambda$, we replace $\bm{u}$ with $\bm{\xi} - \hat{\bm{\xi}}^{(i)}_{h}$ and effectively bound the influence of the deviation $\bm{\xi} - \hat{\bm{\xi}}^{(i)}_{h}$ by scaling it with $\lambda$. Then we interchange the maximization over $\bm{z}_{hij}$ with the minus sign and get an re-expressed upper bound as \eqref{eq:rewrite-interchange}. Finally, we reach \eqref{eq:conjugacy} by applying the definition of conjugacy. Recall the definition of the conjugate function. For a given $\bm{x} \in \mathbb{X}$, we have $[-\ell_{j} + \chi_{\Xi}]^{*}(\bm{v}) = \sup_{\bm{\xi} \in \mathbb{R}^{M}} \left\{\langle \bm{v},\bm{\xi} \rangle+\ell_{j}(\bm{x},\bm{\xi})-\chi_{\Xi}(\bm{\xi})\right\}$. The characteristic function $\chi_{\Xi}$ is defined as $\chi_{\Xi}(\bm{\xi}) = 0$ if $\bm{\xi} \in \Xi$; $\chi_{\Xi}(\bm{\xi}) = \infty$ otherwise. Therefore, the conjugate simplifies to $[-\ell_{j} + \chi_{\Xi}]^{*}(\bm{v}) = \sup_{\bm{\xi} \in \Xi} \left\{\langle \bm{v},\bm{\xi} \rangle+\ell_{j}(\bm{x},\bm{\xi})\right\}$. Subtituting $\bm{v} = -\bm{z}_{hij}$, we have $[-\ell_{j} + \chi_{\Xi}]^{*}(-\bm{z}_{hij}) = \sup_{\xi \in \Xi}\left\{\ell_{j}(\bm{x},\bm{\xi}) - \langle \bm{z}_{hij},\bm{\xi}\rangle\right\}$. The equality from \eqref{eq:conjugacy} to \eqref{eq:substitue-minus-z} holds due to a substitution of $\bm{z}_{hij}$ with $-\bm{z}_{hij}$.

Under Assumption \ref{asp:convexity}, the inequality \eqref{eq:max-min} and \eqref{eq:rewrite-interchange} are reduced to equalities based on Proposition 3.4 in \citep{shapiro2001duality} and Proposition 5.5.4 in \citep{bertsekas2009convex}, leading to the equality between the optimal value of the inner worst-case expectation and \eqref{eq:substitue-minus-z}.

Furthermore, by Theorem 11.23(a) in \citet{rockafellar2009variational}, we have
\begin{align*}
    [-\ell_{j} + \chi_{\Xi}]^{*}(\bm{z}_{hij}) &= \inf_{\bm{\nu}_{hij}}\left([-\ell_{j}]^{*}(\bm{z}_{hij}-\bm{\nu}_{hij}) + [\chi_{\Xi}]^{*}(\bm{\nu}_{hij})\right)\\
    &= \textrm{cl}\left[\inf_{\bm{\nu}_{hij}}\left([-\ell_{j}]^{*}(\bm{z}_{hij}-\bm{\nu}_{hij}) + \sigma_{\Xi}(\bm{\nu}_{hij})\right)\right],
\end{align*}
where $\sigma_{\Xi}$ is the support function of $\Xi$. As the conjugate of $\chi_{\Xi}$, the support function of $\Xi$ is defined via $\sigma_{\Xi}(\bm{\nu}_{hij}):=\sup_{\bm{\xi} \in \Xi}\langle \bm{\nu}_{hij},\bm{\xi}\rangle$.
Since the closure of a function range contains only nonpositive values if and only if the function itself never takes positive values, we conclude that \eqref{eq:substitue-minus-z} is equivalent to \eqref{eq:DRO-finite-convex-program}. This completes the proof.\hfill$\Box$
\endproof

\subsection{Proof of Theorem~\ref{thm:piecewise-affine}}\label{sec:proof-2}
\proof{Proof:}
    Because Assumption \ref{asp:convexity} is satisfied if Assumption \ref{asp:piecewise-affine-objective} holds, the inner worst-case expectation \eqref{eq:worst-case-expectation} is equivalent to \eqref{eq:DRO-finite-convex-program}. From the definition of the conjugacy operator, we have
\begin{align}
    [-\ell_{j}]^{*}(\bm{z}_{hij}) &= \sup_{\bm{\xi}}\langle \bm{z}_{hij},\bm{\xi} \rangle + \langle \bm{a}_{j},\bm{\xi} \rangle + b_{j}\nonumber\\
    &=
    \begin{cases}
    b_{j} & \textrm{if}\quad \bm{z_{hij}} = -\bm{a}_{j},\\\nonumber
    \infty & \textrm{else}.
    \end{cases}
\end{align}
We apply standard duality on the definition of support function
\begin{align}\label{eq:duality-gamma}
\sigma_{\Xi(\bm{\nu}_{hij})} =
\begin{cases}
    \sup_{\bm{\xi}} \quad \langle \bm{\nu}_{hij},\bm{\xi} \rangle \\
    \textrm{s.t.} \quad \bm{C}\bm{\xi} \leq \bm{g}
    \end{cases} = 
    \begin{cases}
    \inf_{\bm{\gamma}_{hij} \geq 0} \quad \langle \bm{\gamma}_{hij},\bm{g} \rangle\\
    \textrm{s.t.} \quad \bm{C}^\mathsf{T}\bm{\gamma}_{hij} = \bm{\nu}_{hij},
\end{cases}
\end{align}
where the equality in \eqref{eq:duality-gamma} follows from strong duality as the uncertainty set is non-empty. This completes the proof.\hfill$\Box$
\endproof

\subsection{Proof of Theorem~\ref{thm:separable-affine}}\label{sec:proof-3}
\proof{Proof:}
We first demonstrate how to reduce the inner worst-case expectation to a finite convex program under all assumptions in this context. The reformulation process up to the point where we introduce epigraphical auxiliary variables follows a similar approach to the proof of Theorem \ref{thm:convex-reduction}, and thus, 
\begin{subequations}\label{eq:interchange-summation-maximization}
\begin{align}
\sup_{\mathbb{Q} \in \mathbb{B}_{\epsilon}(\hat{\mathbb{P}}_{HI})}\mathbb{E}^{\mathbb{Q}}\left[\ell(\bm{x},\bm{\xi})\right] &=
\inf_{\lambda \geq 0} \quad \lambda\epsilon + \sum_{h=1}^{H}\sum_{i=1}^{I} \frac{t^{(I)}_{h}}{I} \sup_{\bm{\xi} \in \Xi} \left\{\ell(\bm{x},\bm{\xi}) - \lambda\|\bm{\xi}-\hat{\bm{\xi}}^{(i)}_{h}\|_{N}\right\},\label{eq:similar-point}\\
&= \inf_{\lambda \geq 0} \quad \lambda\epsilon + \sum_{h=1}^{H}\sum_{i=1}^{I}\sum_{n=1}^{N} \frac{t^{(I)}_{h}}{I} \sup_{\xi_{n} \in \Xi_{n}}\left\{\max_{j \in [J]}\left\{\ell_{nj}(\bm{x},\bm{\xi}_{n})\right\} - \lambda\|\bm{\xi}_{n}-\hat{\bm{\xi}}^{(i)}_{h,n}\|\right\},\label{eq:interchange}
\end{align}
\end{subequations}
where the separability of the overall loss function enable us to interchange the summation and maximization and obtain \eqref{eq:interchange} from \eqref{eq:similar-point}. We then introduce epigraphical auxiliary variables $s_{hin}$, $h \in [H]$, $i \in [I]$, $n \in [N]$:
\begin{subequations}
\begin{align}
    & \quad
    \begin{cases}
        \inf_{\lambda,s_{hin}} \quad \lambda\epsilon + \sum_{h=1}^{H}\sum_{i=1}^{I}\sum_{n=1}^{N}\frac{t^{(I)}_{h}}{I}s_{hin}\\
        \textrm{s.t.} \quad \sup_{\bm{\xi}_{n} \in \Xi_{n}}\left\{\ell_{nj}(\bm{x},\bm{\xi}_{n})-\lambda\|\bm{\xi}_{n} - \hat{\bm{\xi}}^{(i)}_{h,n}\|\right\} \leq s_{hin}, \quad \forall h \in [H], i \in [I], j \in [J], n \in [N]\\
        \quad \quad \lambda \geq 0
    \end{cases}\label{eq:separable-introduce-s}\\
    & =
    \begin{cases}
        \inf_{\lambda,s_{hin},\bm{z}_{hijn}} \quad \lambda\epsilon + \sum_{h=1}^{H}\sum_{i=1}^{I}\sum_{n=1}^{N}\frac{t^{(I)}_{h}}{I}s_{hin}\\
        \textrm{s.t.} \quad \sup_{\bm{\xi}_{n} \in \Xi_{n}}\left\{\ell_{nj}(\bm{x},\bm{\xi}_{n}) - \max_{\|\bm{z}_{hijn}\|_{*} \leq \lambda}\langle \bm{z}_{hijn},\bm{\xi}_{n} - \hat{\bm{\xi}}^{(i)}_{h,n}\rangle\right\} \leq s_{hin},\quad \forall h \in [H], i \in [I], j \in [J], n \in [N]\\
        \quad \quad \lambda \geq 0
    \end{cases}\label{eq:separable-dual-norm-definition}\\
    & \leq
    \begin{cases}
        \inf_{\lambda,s_{hin},\bm{z}_{hijn}} \quad \lambda\epsilon + \sum_{h=1}^{H}\sum_{i=1}^{I}\sum_{n=1}^{N}\frac{t^{(I)}_{h}}{I}s_{hin}\\
        \textrm{s.t.} \quad \sup_{\bm{\xi}_{n} \in \Xi_{n}}\left\{\ell_{nj}(\bm{x},\bm{\xi}_{n})-\langle \bm{z}_{hijn},\bm{\xi}_{n}\rangle\right\}+\langle \bm{z}_{hijn},\hat{\bm{\xi}}^{(i)}_{h,n}\rangle \leq s_{hin}, \quad \forall h \in [H], i \in [I], j \in [J], n \in [N]\\
        \quad \quad \|\bm{z}_{hijn}\|_{*} \leq \lambda, \quad \forall h \in [H], i \in [I], j \in [J], n \in [N]\\
    \end{cases}\label{eq:separable-rewrite-interchange}\\
    & =
    \begin{cases}
        \inf_{\lambda,s_{hin},\bm{z}_{hijn}} \quad \lambda\epsilon + \sum_{h=1}^{H}\sum_{i=1}^{I}\sum_{n=1}^{N}\frac{t^{(I)}_{h}}{I}s_{hin}\\
        \textrm{s.t.} \quad [-\ell_{nj} + \chi_{\Xi_{n}}]^{*}(-\bm{z}_{hijn})+\langle \bm{z}_{hijn},\hat{\bm{\xi}}^{(i)}_{h,n}\rangle \leq s_{hin}, \quad \forall h \in [H], i \in [I], j \in [J], n \in [N],\\
        \quad \quad \|\bm{z}_{hijn}\|_{*} \leq \lambda, \quad \forall h \in [H], i \in [I], j \in [J], n \in [N].\label{eq:separable-conjugacy}\\
    \end{cases}        
\end{align}
\end{subequations}
Here, inequality \eqref{eq:separable-rewrite-interchange} holds in a manner similar to \eqref{eq:rewrite-interchange} and equality is reached given that $\Xi_{n}$ and $\ell_{nj}(\bm{x},\bm{\xi}_{n})$ satisfy the convex assumption in Assumption $\ref{asp:convexity}$. Applying the definition of conjugacy, we obtain \eqref{eq:separable-conjugacy}.

In a similar manner as how we derive \eqref{eq:DRO-piecewise-affine} in Theorem \ref{thm:piecewise-affine}, we have
\begin{align}
    [-\ell_{nj}]^{*}(\bm{z}_{hijn}) &= \sup_{\bm{\xi}_{n}}\langle \bm{z}_{hijn},\bm{\xi}_{n} \rangle + \langle \bm{a}_{nj},\bm{\xi}_{n} \rangle + b_{nj}\nonumber\\
    &=
    \begin{cases}
    b_{nj} & \textrm{if}\quad \bm{z_{hijn}} = -\bm{a}_{nj},\\\nonumber
    \infty & \textrm{else},
    \end{cases}
\end{align}
from the definition of the conjugacy operator.
We apply standard duality on the definition of support function
\begin{align}\label{eq:separable-duality-gamma}
\sigma_{\Xi_{n}(\bm{\nu}_{hijn})} =
\begin{cases}
    \sup_{\bm{\xi}_{n}} \quad \langle \bm{\nu}_{hijn},\bm{\xi}_{n} \rangle \\
    \textrm{s.t.} \quad \bm{C}_{n}\bm{\xi}_{n} \leq \bm{g}_{n}
    \end{cases} = 
    \begin{cases}
    \inf_{\bm{\gamma}_{hijn} \geq 0} \quad \langle \bm{\gamma}_{hijn},\bm{g}_{n} \rangle\\
    \textrm{s.t.} \quad \bm{C}_{n}^\mathsf{T}\bm{\gamma}_{hijn} = \bm{\nu}_{hijn},
\end{cases}
\end{align}
where the equality in \eqref{eq:separable-duality-gamma} follows from strong duality as each uncertainty set $\Xi_{n}$ is non-empty. This completes the proof.
\hfill$\Box$
\endproof
\subsection{Proof of Theorem~\ref{thm:theorem4}}\label{sec:proof-4}
\proof{Proof:}
     Let $T_{i}$ be a time-inhomogeneous Markov chain representing the trust level of the source $1$ with deviation $Y$, defined on the finite state space $S = \{0,\Delta t, 2\Delta t,\ldots,1\}$, $|S| = 1/\Delta t+1$. Note that our assumptions about $\bm{t}^{(0)}$ and $\Delta t$ ensure that $1/\Delta t+1$ is a positive integer. For state $s=0$, the probability of going to state $\Delta t$ in event $i$ is $p^{\text{up}}_{0}(i)$ and the probability of staying in state $0$ is $1-p^{\text{stay}}_{0}(i)$. For state $0<s<1$, the probability of going to state $s+\Delta t$ in event $i$ is $p^{\text{up}}_{s}(i)$, while the probability of going to state $s-\Delta t$ in event $i$ is $p^{\text{down}}_{s}(i) = 1-p^{\text{up}}_{s}(i)-p^{\text{stay}}_{s}(i)$. For state $s=1$, the probability of staying in the current state is $p^{\text{stay}}_{1/\Delta t}(i)$ and the probability of going to state $s - \Delta t$ is $p^{\text{down}}_{s}(i) = 1-p^{\text{stay}}_{1/\Delta t}(i)$.

    If $Y\beta Z$, we know $\mathbb{P}[Y < Z] \geq \beta$, $\beta \geq 0.5$. This is equivalent to $p^{\text{up}}_{s}(i) \geq 0.5$, $\forall s=0,\Delta t,\ldots,1-\Delta t$ and $p^{\text{stay}}_{1}(i) \geq 0$
    for all $i \in [I+1]$. Define the expected change in trust at event $i$ given the current state $T_{i} = s$:
\begin{equation*}
    \mathbb{E}\left[T_{i+1}-T_{i}|T_{i}=s\right] = \Delta t \times (p^{\text{up}}_{s}(i) - p^{\text{down}}_{s}(i)) \geq \Delta t \times (2\beta -1) \geq 0.
\end{equation*}
    We aim to show that $T_{i}$ reaches the maximum trust level $1$ in finite expected time. Since the expected change in trust is at least $t \times (2\beta -1)$, there is a positive drift towards $1$.

    Let $\tau$ be the expected time to reach trust level $1$ from $T_{0}$. Let $V_{i} = 1-T_{i}$ be a nonnegative function representing the ``distance” from maximum trust. Then the expected decrease in $V_{i}$ at each step is:
\begin{equation*}
    \mathbb{E}\left[V_{i+1}-V_{i}|T_{i}=s\right] = -\mathbb{E}\left[T_{i+1}-T_{i}|T_{i}=s\right] \leq -\Delta t \times (2\beta-1).
\end{equation*}
    Thus, $V_{i}$ decreases in expectation by at least $t \times (2\beta-1)$ at each step. Starting from an initial $V_{i}=1-T_{i}$, the expected number of steps to reach $V_{i}=0$ (i.e., $T_{i}=1$) is at most:
\begin{equation*}
    \mathbb{E}[\tau] \leq \frac{V_{0}}{\Delta t \times (2\beta-1)}.
\end{equation*}
Since $V_{0} \leq 1$, $\Delta t >0$, and $2\beta-1>0$, the expected time $\mathbb{E}[\tau]$ is finite.

Once $T_{i}$ reaches $1$, it can only decrease if $Y$ has the maximum error. However, since $p_{1}^{\text{down}}(i) \leq 1-\beta \leq 0.5$, the probability of decreasing is less than or equal to $0.5$. Simiarly, when $T_{i}=1$, the expected change is:
\begin{equation*}
    \mathbb{E}\left[T_{i+1}-T_{i}|T_{i}=1\right] \geq -\Delta t \times (1-\beta).
\end{equation*}
Since $\beta \geq 0.5$, the expected decrease is small, and $T_{i}$ will fluctuate around $1$ within a small interval determined by $\Delta t$ and $\beta$. This completes the proof.\hfill$\Box$
\endproof

\subsection{Proof of Theorem~\ref{thm:theorem5}}\label{sec:proof-5}

\proof{Proof:}
     For source with deviation $Y$, the possible trust changes are: trust increases by $\Delta t$ if $h_{Y} = h^{+}$; trust decreases by $\Delta t$ if $h_{Y} = h^{-}$; trust remains the same otherwise. Therefore, we define $p^{\text{up}}(i) = \mathbb{P}\left[h_{Y} = h^{+}\right]$, $p^{\text{down}}(i) = \mathbb{P}\left[h_{Y} = h^{-}\right]$, $p^{\text{stay}}(i) = 1-p^{\text{up}}(i)-p^{\text{down}}(i)$, for all event $i \in [I+1]$. Our goal is to show that $p^{\text{up}}(i) > p^{\text{down}}(i)$ for all $i \in [I+1]$, ensuring a positive expected change in trust for source $h_{Y}$.

    Since $h_{Y}$ dominates each other source $h \neq h_{Y}$ with probability $\beta$, we have:
\begin{equation*}
    \mathbb{P}[\|\xi^{(i)}_{h_{Y}}\| < \|\xi^{(i)}_{h}\|] \geq \beta, \forall i \in [I+1], \forall h \neq h_{Y}.
\end{equation*}
    Assuming independence between the errors, the probability that $Y$ has the minimum error among all $H$ sources is:
\begin{equation*}
    p^{\text{up}}(i) = \mathbb{P}\left[\cap_{h \neq h_{Y}}(Y<Z_{h})\right] \geq \beta^{H-1}.
\end{equation*}
Similarly, the probability that $Y$ has the maximum error is:
\begin{equation*}
    p^{\text{down}}(i) = \mathbb{P}\left[\cap_{h \neq h_{Y}}(Y>Z_{h})\right] \leq (1-\beta)^{H-1}.
\end{equation*}
Since $\beta \geq 0.5$, we have $p^{\text{up}}(i) \geq 0.5^{H-1}$ and $p^{\text{down}}(i) \leq 0.5^{H-1}$. Thus, $p^{\text{up}}(i) \geq p^{\text{down}}(i)$ and equality holds only when $\beta = 0.5$. Moreover, if $\beta > 0.5$, then $p^{\text{up}}(i) > p^{\text{down}}(i)$. The expected change in trust for $h_{Y}$ at each event is:
\begin{equation*}
    \mathbb{E}\left[T_{i+1}-T_{i}|T_{i}=s\right]=\Delta t \times \delta
\end{equation*}
where $\delta = p^{\text{up}}(i) - p^{\text{down}}(i)$. The expected change is positive, indicating a drift towards maximum trust.

As in the two-source case, thee expected decrease in ``distance” from maximum trust $V_{i} = 1-T_{i}$ is:
\begin{equation*}
    \mathbb{E}\left[V_{i+1}-V_{i}|T_{i}=s\right] = -\mathbb{E}\left[T_{i+1}-T_{i}|T_{i}=s\right] = -\Delta t \times \delta.
\end{equation*}
The expected time $\mathbb{E}[\tau]$ to reach $T_{i}=1$ from $T_{0}$ is bounded by:
\begin{equation*}
    \mathbb{E}[\tau] \leq \frac{V_{0}}{\Delta t \times \delta}.
\end{equation*}

After reaching $T_{i}=1$, the trust in $h_{Y}$ can only decrease if $Y$ has the maximum error. Since $p^{\text{down}}(i))$ is small (especially when $H$ is large and $\beta > 0.5$),  However, since $p_{1}^{\text{down}}(i) \leq 1-\beta \leq 0.5$, the expected decrease is minimal, and $T_{i}$ will fluctuate within a small interval near $1$. This completes the proof.\hfill$\Box$
\endproof

\subsection{Proof of Theorem~\ref{thm:theorem6}}\label{sec:proof-6}

\proof{Proof:}
    Based on Definition \ref{def:definition1}, we have
\begin{align*}
    \mathbb{P}[Y < Z] &= \int \mathbb{P}[Y<a|Z=a]dF_{Z}(a)\\
    &= \int F_{Y|Z}(a)dF_{Z}(a)\\
    &= \mathbb{E}_{Z}[F_{Y|Z}(Z)],
\end{align*}
    where the second equality follows from the continuity of $F_{Y|Z}$. Here, $F_{Y|Z}(a) = \mathbb{P}[Y \leq a | Z = a]$. Given that $Y$ and $Z$ are independent, we immediately have $\mathbb{P}[Y < Z] = \mathbb{E}_{Z}[F_{Y}(Z)]$.
    
    Under the assumption that the smaller the deviation the better, we know $Y$ stochastically dominates $Z$ in the first degree if and only if
    \[F_{Y}(a) \geq F_{Z}(a), \quad\forall a \in [0,+\infty),\]
    \[F_{Y}(a_{0}) > F_{Z}(a_{0}), \quad\exists a_{0} \in [0,+\infty).\]
    And by continuity, there exist $a_{1},a_{2}$ with $a_{1}<a_{0}<a_{2}$, so that
    \[F_{Y}(a) > F_{Z}(a), \quad\forall a \in (a_{1},a_{2}).\]
    Then we have
\begin{align*}
    \mathbb{P}[Y < Z] &= \mathbb{E}_{Z}[F_{Y}(Z)]\\
    &= \int^{a_{1}}_{0}F_{Y}(a)dF_{Z}(a) + \int^{a_{2}}_{a_{1}}F_{Y}(a)dF_{Z}(a) + \int^{+\infty}_{a_{2}}F_{Y}(a)dF_{Z}(a)\\
    &> \int^{a_{1}}_{0}F_{Z}(a)dF_{Z}(a) + \int^{a_{2}}_{a_{1}}F_{Z}(a)dF_{Z}(a) + \int^{+\infty}_{a_{2}}F_{Z}(a)dF_{Z}(a)\\
    &= \int^{+\infty}_{0}F_{Z}(a)dF_{Z}(a) = \mathbb{E}_{Z}F_{Z}(Z).
\end{align*}
    From the probability integral transform, we know that the random variable $F_{Z}(Z)$ has a standard Uniform distribution. Therefore, we have $\mathbb{E}_{Z}F_{Z}(Z) = 0.5$ and $\mathbb{P}[Y < Z] > 0.5$, which implies $\beta = 0.5$. This completes the proof. \hfill$\Box$
\endproof

\subsection{Proof of Theorem~\ref{thm:separable-affine}}\label{sec:proof-7}

\proof{Proof:}
     The Exponential Error Trust Update Algorithm~\ref{alg:exponential-error-trust-update} updates the unnormalized trust weights after event $i$ as
\begin{equation*}
    (t^{\text{new}})^{(i)}_{h} = t^{(i)}_{h}\cdot e^{-\eta \cdot \|\Delta \bm{\xi}^{(i)}_{h}\|}
\end{equation*}
where $\eta >0$ is the update rate. The normalized trust for event $i+1$ is then
\begin{equation*}
    t^{(i+1)}_{h} = \frac{(t^{\text{new}})^{(i)}_{h}}{\sum_{h=1}^{H}(t^{\text{new}})^{(i)}_{h}}.
\end{equation*}
Define the cumulative deviations up to event $i$ for source $h$ as
\begin{equation*}
    S^{(i)}_{h} = \sum_{n=1}^{i}\|\Delta \bm{\xi}^{(n)}_{h}\|.
\end{equation*}
The log-trust of each source can be expressed as:
\begin{equation*}
    \ln (t^{\text{new}})_{h}^{(i)} = \ln t_{h}^{(0)} - \eta S^{(i)}_{h}.
\end{equation*}
We aim to show that $t^{(i)}_{h_{Y}} \rightarrow 1$ as $i \rightarrow \infty$. Consider the difference in log-trust between any other source $h \neq h_{Y}$ and source $h_{Y}$:
\begin{equation*}
     \ln (t^{\text{new}})_{h}^{(i)} -  \ln (t^{\text{new}})_{h_{Y}}^{(i)} = (\ln t_{h}^{(0)} - \ln t_{h_{Y}}^{(0)}) - \eta(S^{(i)}_{h}-S^{(i)}_{h_{Y}}) = - \eta(S^{(i)}_{h}-S^{(i)}_{h_{Y}}).
\end{equation*}

The expected difference in cumulative losses between source $h$ and $h_{Y}$ is:
\begin{equation*}
    \mathbb{E}\left[S^{(i)}_{h}-S^{(i)}_{h_{Y}}\right] = \sum_{n=1}^{i}\left(\mathbb{E}\left[\|\Delta \bm{\xi}^{(n)}_{h}\|-\|\Delta \bm{\xi}^{(n)}_{h_{Y}}\|\right]\right) \geq i\zeta,
\end{equation*}
because for each event $n$, we have $\mathbb{E}\left[\|\Delta \xi^{(i)}_{h}\|\right] - \mathbb{E}\left[\|\Delta \xi^{(i)}_{h_{Y}}\|\right] \leq \zeta$.

Since the deviations $\|\Delta \xi^{(i)}_{h}\|$ are independent and bounded within $[0,L_{\text{max}}]$, we can apply Hoeffding's inequality to bound the probability that the cumulative deviations deviates from its expected value:
\begin{equation*}
    \mathbb{P}\left(S^{(i)}_{h}-S^{(i)}_{h_{Y}} \leq i\zeta - \epsilon\right) \leq \exp\left(-\frac{2\epsilon^{2}}{2iL^{2}_{\text{max}}}\right)
\end{equation*}
This inequality implies that, with high probability, the cumulative loss difference satisfies:
\begin{equation*}
    S^{(i)}_{h}-S^{(i)}_{h_{Y}} \geq i\zeta - o(i). 
\end{equation*}
Then we obtain
\begin{equation*}
   \ln (t^{\text{new}})_{h}^{(i)} -  \ln (t^{\text{new}})_{h_{Y}}^{(i)} = - \eta(S^{(i)}_{h}-S^{(i)}_{h_{Y}}) \leq -\eta(i\zeta - o(i)) = -\eta i\zeta + o(i).
\end{equation*}
Thus, the difference in logarithmic trust between any source $h \neq h_{Y}$ and source $h_{Y}$ becomes increasingly negative as $i$ increases.

The ratio of the unnormalized trust between source $h$ and source $h_{Y}$ is:
\begin{equation*}
    \frac{(t^{\text{new}})_{h}^{(i)}}{(t^{\text{new}})_{h_{Y}}^{(i)}} = e^{\ln (t^{\text{new}})_{h}^{(i)} -  \ln (t^{\text{new}})_{h_{Y}}^{(i)}} \leq e^{-\eta i\zeta + o(i)}.
\end{equation*}
As $i \rightarrow \infty$, the term $e^{-\eta i\zeta + o(i)} \rightarrow 0$ because the exponential of a negative linear term dominates any sublinear terms in $o(i)$. Therefore,
\begin{equation*}
    \lim_{i \rightarrow \infty} \frac{(t^{\text{new}})_{h}^{(i)}}{(t^{\text{new}})_{h_{Y}}^{(i)}} = 0, \quad \forall h \neq h_{Y}.
\end{equation*}
The normalized trust $t^{(i+1)}_{h_{Y}}$ for source $h_{Y}$ would be
\begin{equation*}
    \lim_{i \rightarrow \infty} t^{(i+1)}_{h_{Y}} = \lim_{i \rightarrow \infty} \frac{(t^{\text{new}})_{h_{Y}}^{(i)}}{(t^{\text{new}})_{h_{Y}}^{(i)} + \sum_{h \neq h_{Y}}(t^{\text{new}})_{h}^{(i)}} = 1.
\end{equation*}
This completes the proof.\hfill$\Box$
\endproof
\end{appendices}

\end{document}